\documentclass[11pt,reqno]{amsart}

\usepackage{amsmath,amssymb,amsthm}
\usepackage{bm}
\usepackage{natbib}
\usepackage{graphicx,here,comment}
\usepackage{grffile}
\usepackage{url}

\newtheorem{theorem}{Theorem}

\newtheorem{assumption}{Assumption}
\theoremstyle{definition}

\newtheorem{remark}{Remark}

\newcommand{\R}{\mathbb{R}}
\newcommand{\mF}{\mathcal{F}}

\newcommand{\mA}{\mathcal{A}}
\newcommand{\mB}{\mathcal{B}}

\newcommand{\Ep}{\mathrm{E}}
\renewcommand{\Pr}{\mathrm{P}}
\newcommand{\Prr}{\mathrm{P}}
\newcommand{\mE}{\mathcal{E}}
\renewcommand{\hat}{\widehat}
\renewcommand{\tilde}{\widetilde}

\DeclareMathOperator{\Var}{Var}
\DeclareMathOperator{\Cov}{Cov}

\usepackage{color}

\newcommand{\bc}{\color{blue}}

\begin{document}

\title[Functional linear regression with functional responses]{PCA-based estimation for functional linear regression with functional responses}
\thanks{M. Imaizumi is supported by Grant-in-Aid for JSPS Research Fellow (15J10206) from the JSPS. K. Kato is supported by Grant-in-Aid for Scientific Research (C) (15K03392) from the JSPS}

\author[M. Imaizumi]{Masaaki Imaizumi}
\author[K. Kato]{Kengo Kato}

\date{First version: August 26, 2016. This version: \today}

\address[M. Imaizumi]{
Graduate School of Economics, University of Tokyo,
7-3-1 Hongo, Bunkyo-ku, Tokyo 113-0033, Japan.
}
\email{imaizumi-masaaki@g.ecc.u-tokyo.ac.jp}

\address[K. Kato]{
Graduate School of Economics, University of Tokyo,
7-3-1 Hongo, Bunkyo-ku, Tokyo 113-0033, Japan.
}
\email{kkato@e.u-tokyo.ac.jp}

\begin{abstract}
This paper studies a regression model where both predictor and response variables are random functions. 
We consider a functional linear model where the conditional mean of the response variable at each time point is given by a linear functional of the predictor variable. 
In this paper, we are interested in estimation of the integral kernel $b(s,t)$ of the conditional expectation operator, where $s$ is an output variable while $t$ is a variable that interacts with the predictor variable. This problem is an ill-posed inverse problem, and we consider two estimators based on the functional principal component analysis (PCA). We  show that under suitable regularity conditions, an estimator based on the single truncation  attains the convergence rate for the integrated squared error that is characterized by  smoothness of the function $b (s,t)$ in $t$ together with the decay rate of the eigenvalues of the covariance operator, but the rate does not depend on smoothness of $b(s,t)$ in $s$. This rate is shown to be minimax optimal, and consequently smoothness of $b(s,t)$ in $s$ does not affect difficulty of estimating $b$.
We also consider an alternative estimator based on the double truncation, and provide conditions under which the alternative estimator attains the optimal rate. 
We conduct simulations to verify the performance of PCA-based estimators in the finite sample. Finally, we apply our estimators to investigate the relation between the lifetime pattern of working hours and total income, and the relation between the electricity spot price and the wind power infeed.
\end{abstract}



\maketitle

\section{Introduction}

This paper studies a regression model where both predictor and response variables are random functions. Let $X,Y$ be $L^{2}(I)$-valued random variables with $I=[0,1]$, and  consider a regression model of the form 
\begin{equation}
\Ep(Y \mid X)(s) = \Ep\{Y(s)\} + \int_{I} b(s,t) [ X(t) - \Ep\{X(t)\} ] dt. \label{eq: model0}
\end{equation}
See Section \ref{sec: setup} for the precise description of the setup. 
The focus of this paper is on estimation of the bivariate function $b(s,t)$, which is an ill-posed inverse problem (see Remark \ref{rem: ill-posed inverse problem} in Section \ref{sec: setup}).

Data collected on dense grids can be typically regarded as realizations of a random function (stochastic process), and such data are called \textit{functional data}. 
Statistical methodology dealing with functional data is called \textit{functional data analysis} and has a large number of fruitful applications \citep[see][]{RaSi05}. 
For example, the functional linear model (\ref{eq: model0}) with functional predictor and response variables can be used to investigate how a complete daily temperature profile over one year influences a daily precipitation at each day \citep[][Chapter 16]{RaSi05}.

In this paper, we consider estimators for the function $b$  based on the functional principal component analysis (PCA), which is one of standard techniques used in functional data analysis. Applying basis expansions of $X$ and $b$ using  the eigenfunction system $\{ \phi_{k} \}_{k=1}^{\infty}$ for the covariance operator of $X$, we can expand $X$ and $b$ as $X(t)= \Ep\{X(t)\}+ \sum_{k} \xi_{k} \phi_{k}(t)$ and $b (s,t)= \sum_{j,k} b_{j,k} \phi_{j}(s)\phi_{k}(t)$, where we measure smoothness of $b$ via how fast $|b_{j,k}|$ decays as $j \to \infty$ or $k \to \infty$.  We consider two methods to estimate $b$ based on different characterizations of $b$. The first method uses the fact that $\Ep\{ \xi_{k}Y(s) \} = \Ep(\xi_{k}^{2}) \sum_{j} b_{j,k} \phi_{j}(s)$. This method is based on truncation of the series expansion $b(s,t) = \sum_{k} [ \Ep \{\xi_{k}Y(s)\}/\Ep(\xi_{k}^{2}) ] \phi_{k}(t)$  by a finite series $\sum_{k=1}^{m_{n}}$ with $m_{n} \to \infty$ as $n \to \infty$ (which we call the single truncation in comparison with the second method below), and replace $\Ep \{\xi_{k}Y(\cdot)\}$, $\Ep(\xi_{k}^{2})$, and $\phi_{k}$ by their estimators.
This estimator was considered by \cite{CrMa13}. The second method uses the expansion of $Y$ as $Y(s) = \sum_{j} \eta_{j} \phi_{j}(s)$. This alternative method is based on truncation of the double series expansion $b(s,t) = \sum_{j,k} \{ \Ep(\eta_{j}\xi_{k})/\Ep(\xi_{k}^{2}) \}  \phi_{j}(s)\phi_{k}(t)$ by a finite series $\sum_{j=1}^{m_{n,1}} \sum_{k=1}^{m_{n,2}}$ with $m_{n,1} \to \infty$ and $m_{n,2} \to \infty$ as $n \to \infty$  (which we call the double truncation), and replace $\Ep(\eta_{j}\xi_{k}), \Ep(\xi_{k}^{2})$, and $\phi_{j}$ by their estimators.

\cite{CrMa13} consider our first estimator, but the focus in \cite{CrMa13} is on prediction, and not on estimation of the function $b$ \textit{per se}. These two problems are substantially different, and they do not derive sharp rates of convergence for their estimator of $b$ itself.
\cite{PaQi12} and \cite{HoKi15} analyze the estimator of \cite{CrMa13} for $b$ with dependent functional data, but they only prove consistency of the estimator.  \cite{YaMuWa05} consider a PCA-based estimator similar to our second estimator, but do not explicitly derive rates of convergence for their estimator. 

The object of this paper is to study rates of convergence for estimation of $b$. 
First, we show that under suitable regularity conditions, the estimator based on the single truncation (that is, the estimator of \cite{CrMa13}) attains the convergence rate for the integrated squared error that is characterized by  smoothness of the function $b (s,t)$ in $t$ together with the the decay rate of the eigenvalues of the covariance operator, but the rate does not depend on smoothness of $b(s,t)$ in $s$. This rate is shown to be minimax optimal. 
This means that smoothness of $b(s,t)$ in $s$ does not affect difficulty of estimating $b$, which is in sharp contrast with nonparametric estimation of a bivariate regression function.
Next, we analyze the second estimator based on the double truncation, and provide conditions under which it attains the optimal rate. 
We point out that some restrictions on smoothness levels for $b(s,t)$ in $s$ and $t$ are required for the second estimator to achieve the optimal rate. We include the analysis of the second estimator since in applications, the double truncation typically leads to an estimate more interpretable than the single truncation, although from a theoretical point of view, the single truncation is enough for the purpose of estimating $b$; see Remark \ref{rem: double truncation} ahead and the discussion in Chapter 16 of \cite{RaSi05}. We also conduct simulations to verify the performance of the estimators in the finite sample. Finally, we apply our estimators to investigate two topics: the relation between the lifetime pattern of working hours and total income using the data from National Longitudinal Survey of Youth conducted by \cite{NLSY12}, and the relation between the hourly electricity spot prices and the amount of the wind power infeed using the data from EEX Transparency Platform introduced in \cite{Li13}.

The literature on functional data analysis is now quite broad. We refer to \cite{Bo00}, \cite{RaSi05}, and \cite{HsEu15} as general references on functional data analysis. 
One of the main focuses in the previous literature on functional data analysis is a functional linear model with a scalar response variable. See \cite{CaFeSa99,CaFeSa03}, \cite{CaHa06}, \cite{HaHo07}, \cite{LiHs07}, \cite{CrKnSa09}, \cite{JaWaZh09}, \cite{YuCa10}, \cite{CaJo10}, \cite{CaYu12}, \cite{DeHa12}, and \cite{CoJo12}. In particular, \cite{HaHo07} consider a PCA-based estimator and an  estimator based on Tikhonov regularization for the slope function, and provide conditions under which those estimators attain minimax rates of convergence for the integrated squared error. 

The analysis of functional responses was first considered by \cite{RaDa91}. \cite{ChMuWa04} consider a regression model where a predictor variable is finite-dimensional while a response variable is a random function.  
Functional linear models with functional predictor and response variables are considered in 
\cite{CuFeFr02}, \cite{YaMuWa05}, \cite{HeMuWaYa10}, \cite{CrMa13}, \cite{Li15}, \cite{HoKi15}, and \cite{BeCaFl15}. 
\cite{CuFeFr02} work with fixed designs, which is a different setting than ours, and prove consistency of a series estimator of the integral operator with kernel $b$ for  the operator norm.
We already referred to \cite{YaMuWa05}, \cite{CrMa13}, and \cite{HoKi15}. \cite{HeMuWaYa10} propose an estimator of $b$ based on the functional canonical correlation analysis, but do not study its asymptotic properties.
\cite{Li15} considers prediction for functional linear regression with functional responses based on a reproducing kernel Hilbert space approach, which is a topic substantially different from ours.
The recent preprint by \cite{BeCaFl15} studies a Tikhonov regularization estimation for $b$ and establishes rates of convergence for their estimator; the estimator and the assumptions in \cite{BeCaFl15} are substantially different from ours and so their results are not directly comparable to ours. 

Importantly, none of these papers derives optimal rates of convergence for estimation of $b$;  the present paper fills this important void and thereby contributes to advancing the understanding of functional data analysis. 
From a technical point of view, the proofs of the main theorems (Theorems \ref{thm: main1} and \ref{thm: main2}) build upon the techniques developed in \cite{HaHo07}. 
However, since we are estimating a bivariate function with two different levels of smoothness rather than a univariate function in the scalar response case, the proofs require a chain of delicate calculations.
Furthermore, to establish minimax lower bounds for estimating $b$, we have to construct a suitable sequence of conditional distributions of $Y$ given $X$, and since $Y$ takes values in $L^{2}(I)$, we have to construct a sequence of distributions on $L^{2}(I)$, which is a significant difference from \cite{HaHo07}. To this end, we employ the theory of Gaussian measures on Banach spaces \citep[cf.][Chapter VIII]{St11}. 

In this paper, we use basic results on functional analysis. We refer to \cite{ReSi80}  as a general reference on functional analysis. 
\cite{Bo00} and \cite{HsEu15} cover results on functional analysis useful for functional data analysis. 
For mathematical background on linear inverse problems, we refer to \cite{Kr99}. 

The rest of the paper is organized as follows. 
In Section \ref{sec: setup}, we formally describe the setup and estimators. 
In Section \ref{sec: main results}, we present the main results on rates of convergence of the PCA-based estimators for the coefficient function. 
In Section \ref{sec: simulation}, we present simulation results to verify performance of the PCA-based estimates in the finite sample. In Section \ref{sec: real data analysis}, we present applications of our estimators to two real data examples. 
All the proofs are deferred to Appendix. 

\subsection{Notation}
We use the following notation. For any measurable functions $f: I \to \R$ and $R: I^{2} \to \R$, let
$\| f \| = \left\{\int_{I} f^{2}(t) dt\right\}^{1/2}$ and $||| R |||=\left\{\iint_{I^2} R^{2}(s,t) dsdt\right\}^{1/2}$.
For any functions $f,g: I \to \R$, define $f \otimes g: I^{2} \to \R$ by $(f \otimes g)(s,t) = f(s)g(t)$ for $s,t \in I$. 
Let $\mathcal{L}^{2}(I) = \{ f: I \to \R : f \ \text{is measurable}, \ \| f \| < \infty \}$,
and define the equivalence relation $\sim$ for real-valued functions $f,g$ defined on $I$ by $f \sim g \Leftrightarrow f=g$ almost everywhere. Define $L^{2}(I)$ by the quotient space $L^{2}(I) = \mathcal{L}^{2}(I) \slash \sim$ equipped with the inner product $\langle  f^{\sim},g^{\sim}\rangle = \int_{I} f(t)g(t)dt$ for $f,g \in \mathcal{L}^{2}(I)$ where $f^{\sim} = \{ h \in \mathcal{L}^{2}(I) : h \sim f \}$; the space $L^{2}(I)$ is a separable Hilbert space, and as usual, we identify any element in $\mathcal{L}^{2}(I)$ as an element of $L^{2}(I)$. 
Define $L^{2}(I^{2})$ analogously. We also identify any real-valued function  $f$ defined almost everywhere on $I$ (or $I^{2}$) as a function defined everywhere on $I$ (or $I^{2}$) by setting $f(t)=0$ for any point $t$ at which $f$ is not defined. 
For any positive sequences $a_{n},c_{n}$, we write $a_{n} \sim c_{n}$ if $a_{n}/c_{n}$ is bounded and bounded away from zero. 
In what follows, let $(\Omega, \mA, \Prr)$ denote an underlying probability space. 

\section{Setup and estimators}
\label{sec: setup}

Suppose that we observe a pair of random functions $(X,Y)$ indexed by $I =[0,1]$ where $X=\{ X(t) : t \in I \}$ and $Y=\{ Y(t) : t \in I \}$ are predictor and response variables, respectively.  
We assume that $X$ and $Y$ are $L^{2}(I)$-valued random variables such that $\Ep (\| X \|^{2}) < \infty$ and $\Ep (\| Y \|^{2}) < \infty$ (recall that a measurable stochastic process with paths in $L^{2}(I)$ almost surely induces an $L^{2}(I)$-valued random variable, and vice versa; see \cite{Ra72} or \cite{By77}).
We consider a functional linear regression model
\begin{equation}
\Ep( Y \mid X )(s) = \Ep \{Y(s)\}  + \int_{I}  b(s,t) [ X(t) - \Ep \{X(t)\} ]  dt, \label{eq: model}
\end{equation}
where $\Ep( Y \mid X )$ is the conditional expectation of $Y$ as an $L^{2}(I)$-valued random variable conditionally on the $\sigma$-field generated by $X$ (which is well-defined since $\Ep( \| Y \| ) < \infty$, and $\Ep( Y \mid X )$ itself is an $L^{2}(I)$-valued random variable; see Chapter 5 in \cite{St11}), and $(s,t) \mapsto b(s,t)$ is the coefficient function assumed to be in $L^{2}(I^{2})$, that is, $||| b |||^{2} =  \iint_{I^2} b^{2}(s,t) ds dt < \infty$.
The equality in (\ref{eq: model}) should be understood as an equality as $L^{2}(I)$-valued random variables. 

The goal of this paper is estimation of the function $(s,t) \mapsto b(s,t)$, and to this end we shall employ the functional principal component analysis (PCA).
Consider the covariance function 
\[
K(s,t) = \Cov \{X(s),X(t)\}, \ s,t \in I.
\]
The assumption that $\Ep (\|X\|^{2}) < \infty$ ensures that $K \in L^{2}(I^{2})$. 
In addition, we assume that the integral operator from $L^{2}(I)$ into itself with kernel $K$, namely the covariance operator of $X$, is injective (which is equivalent to the condition that $\Var (\langle f,X \rangle) > 0$ for all $f \in L^{2}(I)$ with $\| f \| =1$). The covariance operator is self-adjoint and positive definite. The Hilbert-Schmidt theorem \citep[see][Theorem VI.16]{ReSi80} then ensures that $K$ admits the spectral expansion 
\[
K(s,t) = \sum_{k=1}^{\infty} \kappa_{k} \phi_{k}(s) \phi_{k}(t)
\]
in $L^{2}(I^{2})$, where $\kappa_{1} \geq \kappa_{2} \geq \cdots > 0$ are a non-increasing sequence of eigenvalues tending to zero and $\{ \phi_{k} \}_{k=1}^{\infty}$ is an orthonormal basis of $L^{2}(I)$ consisting of  eigenfunctions of the integral operator, namely, $\int_{I} K(s,t) \phi_{k}(t) dt = \kappa_{k} \phi_{k} (s)$
for all $k \geq 1$.
We will later assume that there are no ties in $\kappa_{j}$'s, that is, $\kappa_{1} > \kappa_{2} > \cdots > 0$. 
Since $\{ \phi_{k} \}_{k=1}^{\infty}$ is an orthonormal basis of $L^{2}(I)$, we have the following expansion in $L^{2}(I)$:
\[
X(t) = \Ep\{ X(t) \} + \sum_{k=1}^{\infty} \xi_{k} \phi_{k}(t),
\]
where each $\xi_{k}$ is defined by 
\[
\xi_{k} = \int_{I} [ X(t) - \Ep\{X(t)\} ] \phi_{k}(t) dt.
\]
By Parseval's identity and Fubini's theorem, $\sum_{k=1}^{\infty} \Ep(\xi_{k}^{2} ) = \int_{I} \Var \{X(t)\} dt < \infty$ and
\begin{equation}
\Ep( \xi_{k} \xi_{\ell} ) = \iint_{I^{2}} K(s,t) \phi_{k}(s) \phi_{\ell}(t) ds dt = 
\begin{cases}
\kappa_{k} & \text{if} \ k=\ell \\
0 & \text{if} \ k \neq \ell
\end{cases}
.
\label{eq: PC score}
\end{equation}

Furthermore, since $\{ \phi_{j} \otimes \phi_{k} \}_{j,k=1}^{\infty}$ is an orthonormal basis of  $L^{2}(I^{2})$, we have 
\[
b(s,t) = \sum_{j,k=1}^{\infty} b_{j,k} \phi_{j}(s)\phi_{k}(t)
\]
in $L^{2}(I^{2})$ with $b_{j,k} = \iint_{I^{2}} b(s,t) \phi_{j} (s) \phi_{k} (t) ds dt$.
This yields that 
\[
\int_{I} b(s,t) [ X(t) - \Ep\{X(t)\} ] dt = \sum_{j=1}^{\infty}  \left( \sum_{k=1}^{\infty} b_{j,k} \xi_{k} \right)\phi_{j}(s).
\]

Now, because of (\ref{eq: PC score}) and since the expansion of $X$ holds in $L^{2}(I \times \Omega, dt \otimes d\Prr)$ too (that is, $\Ep[ \| X - \Ep\{X(\cdot)\} - \sum_{k=1}^{N} \xi_{k} \phi_{k}\|^{2}] = \sum_{k=N+1}^{\infty} \Ep(\xi_{k}^{2}) \to 0$ as $N \to \infty$), we have that $\Ep\{ \xi_{k} Y(s) \} = \kappa_{k} \sum_{j=1}^{\infty} b_{j,k} \phi_{j}(s)$, 
where the equality holds in $L^{2}(I)$, and therefore we obtain the following characterization of $b$:
\begin{equation}
b(s,t) = \sum_{k=1}^{\infty} \frac{\Ep\{\xi_{k}Y(s)\}}{\kappa_{k}} \phi_{k}(t). \label{eq: first characterization}
\end{equation}
This characterization leads to a method to estimate $b$. 

Let $(X_{1},Y_{1}),\dots,(X_{n},Y_{n})$ be independent copies of $(X,Y)$ as ($L^{2}(I) \times L^{2}(I)$)-valued random variables. 
We estimate $K$ by the empirical covariance function $\hat{K}$ defined as
\[
\hat{K} (s,t)= \frac{1}{n} \sum_{i=1}^{n} \{ X_{i}(s)- \overline{X}(s) \} \{ X_{i}(t)- \overline{X}(t) \}, \ s,t \in I,
\]
where $\overline{X} = n^{-1}\sum_{i=1}^{n}X_{i}$. Let 
\begin{equation}
\hat{K}(s,t) = \sum_{k=1}^{\infty} \hat{\kappa}_{k} \hat{\phi}_{k}(s) \hat{\phi}_{k}(t)
\label{eq: expansion of Khat} 
\end{equation}
be the spectral expansion of $\hat{K}$ in $L^{2}(I^{2})$, where $\hat{\kappa}_{1} \geq \hat{\kappa}_{2} \geq \cdots  \geq 0$ are a non-increasing sequence of eigenvalues tending to zero and $\{ \hat{\phi}_{k} \}_{k=1}^{\infty}$ is an orthonormal basis of $L^{2}(I)$ consisting eigenfunctions of the integral operator with kernel $\hat{K}$, namely, $\int_{I} \hat{K}(s,t) \hat{\phi}_{k}(t) dt = \hat{\kappa}_{k} \hat{\phi}_{k}(s)$
for all $k \geq 1$.
The spectral expansion in (\ref{eq: expansion of Khat}) is possible since the integral operator with kernel $\hat{K}$ is of finite rank (at most $(n-1)$), and so in addition to an orthonormal system of $L^{2}(I)$ consisting of eigenfunctions corresponding to the positive eigenvalues, we can add functions so that the augmented system of functions $\{ \hat{\phi}_{k} \}_{k=1}^{\infty}$ becomes an orthonormal basis of $L^{2}(I)$. 

Furthermore, let 
\[
\hat{\xi}_{i,k} = \int_{I} \{ X_{i}(t)-\overline{X}(t) \} \hat{\phi}_{k}(t)dt.
\]
Using the characterization in (\ref{eq: first characterization}), we consider the following estimator based on the single truncation:
\begin{equation}
\hat{b} (s,t) = \sum_{k=1}^{m_{n}} \frac{n^{-1}\sum_{i=1}^{n} \hat{\xi}_{i,k}Y_{i}(s)}{\hat{\kappa}_{k}} \hat{\phi}_{k}(t), \label{eq: first estimator}
\end{equation}
where $m_{n} \to \infty$ as $n \to \infty$. 
This estimator was considered in \cite{CrMa13}. 

We also consider an alternative estimator based on truncating the double series, namely, the double truncation. Let $\mE = Y-\Ep(Y \mid X)$, and consider the expansions $Y(s) = \sum_{j=1}^{\infty} \eta_{j} \phi_{j}(s)$ and $\mE (s) = \sum_{j=1}^{\infty} \varepsilon_{j} \phi_{j}(s)$
in $L^{2}(I)$. Now, since $\int_{I} (\int_{I}  b(s,t) [ X(t) - \Ep\{ X(t) \} ] dt ) \phi_{j} (s) ds= \sum_{k=1}^{\infty} b_{j,k} \xi_{k}$ for each $j \geq 1$,
we have that
\begin{equation}
\eta_{j} = a_{j} + \sum_{k=1}^{\infty} b_{j,k} \xi_{k} + \varepsilon_{j}, \ j \geq 1, \label{eq: model2}
\end{equation}
where $a_{j} = \Ep (\eta_{j}) = \int_{I} \Ep\{Y(s)\} \phi_{j}(s)ds$ for $j \geq 1$.  Therefore, we have $\Ep( \eta_{j} \xi_{k} ) = b_{j,k} \Ep( \xi_{k}^{2} ) = \kappa_{k} b_{j,k}$,
namely, 
\begin{equation}
b_{j,k} = \Ep(\eta_{j}\xi_{k})/\kappa_{k}. \label{eq: representation of b}
\end{equation}
Based on this characterization, we consider the following alternative estimator:
\begin{equation}
\tilde{b}(s,t) = \sum_{j=1}^{m_{n,1}} \sum_{k=1}^{m_{n,2}} \tilde{b}_{j,k} \hat{\phi}_{j}(s) \hat{\phi}_{k}(t), \label{eq: second estimator}
\end{equation}
where $m_{n,1} \to \infty$ and $m_{n,2} \to \infty$ as $n \to \infty$, and each $\tilde{b}_{j,k}$ is defined by $\tilde{b}_{j,k} = n^{-1} \sum_{i=1}^{n} \hat{\eta}_{i,j} \hat{\xi}_{i,k}/\hat{\kappa}_{k}$ with $\hat{\eta}_{i,j} = \int_{I} Y_{i}(s) \hat{\phi}_{j}(s) ds$.

In the next section, we will derive rates of convergence of the estimators $\hat{b}$ and $\tilde{b}$ for the integrated squared error. 

\begin{remark}[Motivation of the double truncation]
\label{rem: double truncation}
It will turn out in the next section that $\hat{b}$ with properly chosen $m_{n}$ is rate optimal, and from a theoretical point of view, the single truncation is enough for the purpose of estimating $b$. 
However, in practice, the double truncation would be a preferred option since, compared with the single truncation, the double truncation typically results in an estimate of $b(s,t)$ more regular in $s$ and thereby yielding a more interpretable estimate. See the discussion in Chapter 16 of \cite{RaSi05} and the real data analysis in Section \ref{sec: real data analysis}.
Hence the analysis of our second estimator is of some importance. 
\end{remark}

\begin{remark}[Ill-posedness of estimation of $b$]
\label{rem: ill-posed inverse problem}
The problem of estimating $b$ can be regarded as a problem of estimating an unknown operator in the operator equation, and therefore is an ill-posed inverse problem. 
For any $R \in L^{2}(I^{2})$, let $T_{R}: L^{2}(I) \to L^{2}(I)$ denote the integral operator with kernel $R$, i.e, 
\[
(T_{R} h)(s) = \int_{I} R(s,t) h(t) dt, \ h \in L^{2}(I).
\]
The adjoint operator $T_{R}^{*}$ of $T_{R}$ is also an integral operator and  of the form 
\[
(T_{R}^{*}h)(t) = \int_{I} R(s,t) h(s) ds, \ h \in L^{2}(I).
\]
Now, let $C_{XY}(s,t) = \Cov \{X(s),Y(t)\}, s,t \in I$. Then, using the symmetry of $K$,  we have that for any $h \in L^{2}(I)$, 
\[
(T_{C_{XY}} h)(t) = \iint_{I^{2}} K(t,u) b(s,u) h(s) ds du = (T_{K} T_{b}^{*}h)(t),
\]
that is, $T_{C_{XY}} =T_{K}T_{b}^{*}$.
Since we are assuming that $T_{K}$ is injective, we have that $T_{b}^{*} = T_{K}^{-1}T_{C_{XY}}$. Both $\Cov \{X(s),Y(t)\}$ and $K$ can be directly estimated from the data. However, since $T_{K}$ is a compact operator \citep[][Theorems VI.22 and VI.23]{ReSi80}, $T_{K}^{-1}$ is necessarily unbounded \citep[][p.23]{Kr99}, and therefore the problem of recovering $T_{b}^{*}$ is ill-posed \citep[][Section 15.1]{Kr99}. 
In fact, consider $C_{XY}^{N} = C_{XY} + \kappa_{N} \phi_{N} \otimes \phi_{1}$, which converges to $C_{XY}$ in $L^{2}(I^{2})$ as $N \to \infty$ (that is, $T_{C^{N}_{XY}}$ converges to $T_{C_{XY}}$ in the Hilbert-Schmidt norm). 
It is seen that $b^{N} = b + \phi_{1} \otimes \phi_{N}$ satisfies $T_{b^{N}}^{*} = T_{K}^{-1} T_{C_{XY}^{N}}$, but $||| b^{N} - b ||| = ||| \phi_{1} \otimes \phi_{N} ||| = 1$.

\end{remark}

\section{Main results}
\label{sec: main results}

\subsection{Rates of convergence}

In this subsection, we derive rates of convergence of the estimators $\hat{b}$ and $\tilde{b}$ defined in (\ref{eq: first estimator}) and (\ref{eq: second estimator}), respectively. To this end, we make the following assumption. Recall that $\mathcal{E}=Y-\Ep(Y \mid X)$. 

\begin{assumption}
\label{as: 1}
There exist constants $\alpha >1, \beta > \alpha/2+1, \gamma > 1/2$, and $C_{1} > 1$ such that
\begin{align}
&\Ep (\| Y \|^{2}) < \infty, \ \Ep( \| X \|^{2}) < \infty, \ \Ep (\| \mathcal{E} \|^{2} \mid X)\leq C_{1} \ \text{almost surely}, \label{eq: moment condition1} \\
&\Ep( \xi_{k}^{4} ) \leq C_{1} \kappa_{k}^{2}, \ \text{for all} \ k \geq 1, \label{eq: moment condition2} \\
&\kappa_{k} \leq C_{1} k^{-\alpha}, \ \kappa_{k} - \kappa_{k+1} \geq C_{1}^{-1} k^{-\alpha-1}, \ \text{for all} \ k \geq 1, \label{eq: eigenvalue condition} \\
&| b_{j,k} | \leq C_{1}j^{-\gamma} k^{-\beta}, \ \text{for all} \ j,k \geq 1.  \label{eq: condition on b}
\end{align}
\end{assumption}

Some comments on Assumption \ref{as: 1} are in order. The first row (\ref{eq: moment condition1}) is a standard moment condition. 
The second (\ref{eq: moment condition2}) and third rows (\ref{eq: eigenvalue condition}) are adapted from \cite{HaHo07}. 
Condition (\ref{eq: moment condition2}) is standard in the literature on functional linear models. Concretely, Condition (\ref{eq: moment condition2}) is automatically satisfied if $X$ is Gaussian, since in that case $\xi_{k}$ are Gaussian. 
In Condition (\ref{eq: eigenvalue condition}), as in \cite{CaHa06} and \cite{HaHo07}, we require that the eigenvalues $\{ \kappa_{k} \}_{k=1}^{\infty}$ are ``well-separated'', namely, $\kappa_{k} - \kappa_{k+1} \geq C_{1}^{-1}k^{-\alpha-1}$ for all $k \geq 1$. This condition is used to ensure sufficient estimation accuracy of the empirical eigenfunctions $\hat{\phi}_{k}$. 
This condition also ensures that, since $\kappa_{k} \to 0$ as $k \to \infty$, $\kappa_{k} = \sum_{j=k}^{\infty} (\kappa_{j} - \kappa_{j+1}) \geq C_{1}^{-1} \sum_{j=k}^{\infty} j^{-\alpha-1} \geq k^{-\alpha}/(C_{1}\alpha)$.
So $\kappa_{k} \sim k^{-\alpha}$ as $k \to \infty$. The value of $\alpha$ measures ``ill-posedness'' of the estimation problem, so that the larger $\alpha$ is, the more difficult estimation of $b$ will be. 
For given constants $\alpha > 1$ and $C_{1} > 1$, the class of distributions of $X$ verifying (\ref{eq: moment condition1})--(\ref{eq: eigenvalue condition}) is rich enough, and a superset of the subclass
\begin{align*}
\Bigg \{ \Pr \circ X^{-1} &: X = \sum_{k} \sqrt{\kappa_{k}} U_{k} \phi_{k}, \ \{ \phi_{k} \} \ \text{is an orthonormal basis of} \ L^{2}(I), \\
&\quad  \{ U_{k} \} \sim WN(0,1), \ \Ep[ U_{k}^{4} ] \leq C_{1} \ \forall k \geq 1, \\
&\quad \kappa_{k} \leq C_{1} k^{-\alpha}, \ \kappa_{k} - \kappa_{k+1} \geq C_{1}^{-1} k^{-\alpha-1} \ \forall k \geq 1 \Bigg \},
\end{align*}
where $\{ U_{k} \} \sim WN(0,1)$ means that $\{ U_{k} \}$ is a white noise process (i.e., an uncorrelated sequence of random variables) with mean zero and unite variance. 

The last condition (\ref{eq: condition on b}) is a smoothness condition on $b$, where the smoothness is measured through the eigenfunction system $\{ \phi_{k} \}_{k=1}^{\infty}$, which is natural in our setting. Since $b(s,t)$ is a bivariate function, however, there are potentially a number of variations on how $b_{j,k}$ decays as $j \to \infty$ or $k \to \infty$.
We focus on a simple case where $|b_{j,k}|$ decays like $j^{-\gamma}k^{-\beta}$ as $j \to \infty$ or $k \to \infty$, and $\gamma$ measures smoothness of $b (s,t)$ in $s$ while $\beta$ measures smoothness of $b(s,t)$ in $t$.   We also require that $\beta > \alpha/2+1$ for a technical reason; see the discussion after Theorem \ref{thm: main1}. 

The following theorem establishes rates of convergence for $\hat{b}$. 

\begin{theorem}
\label{thm: main1}
Consider the estimator $\hat{b}$ defined in (\ref{eq: first estimator}). Suppose that Assumption \ref{as: 1} is satisfied. 
Choose $m_{n}$ in such a way that $m_{n} \to \infty$ and $m_{n} = o\{n^{1/(2\alpha+2)}\}$. 
Then 
\begin{equation}
||| \hat{b} - b |||^{2} = O_{\Prr} \{n^{-1}m_{n}^{\alpha+1} + m_{n}^{-2\beta+1}\}. \label{eq: rate}
\end{equation}
Therefore, by choosing $m_{n} \sim n^{1/(\alpha+2\beta)}$, we have 
\[
||| \hat{b} - b |||^{2} = O_{\Prr} \{n^{-(2\beta-1)/(\alpha+2\beta)}\}.
\]
\end{theorem}

\begin{remark}
It is not difficult to verify from the proof of Theorem \ref{thm: main1} that the results of the theorem hold uniformly over  a class of distributions $\mF (\alpha,\beta,\gamma,C_{1})$ of $(X,Y)$ that verify (\ref{eq: model}) and (\ref{eq: moment condition1})--(\ref{eq: condition on b}) for given constants $\alpha>1,\beta>\alpha/2+1,\gamma>1/2$, and $C_{1} > 1$. In particular, by choosing  $m_{n}\sim n^{1/(\alpha+2\beta)}$, we have 
\[
\lim_{D \to \infty} \limsup_{n \to \infty} \sup_{F \in \mathcal{F}(\alpha,\beta,\gamma,C_{1})} \Prr_{F} \left\{ ||| \hat{b} - b |||^{2} > D n^{-(2\beta-1)/(\alpha+2\beta)} \right\} = 0,
\]
where $\Prr_{F}$ denotes the probability under $F$. We will show in Theorem \ref{thm: lower bound}  that the rate $n^{-(2\beta-1)/(\alpha+2\beta)}$ is minimax optimal. 
\end{remark}

The requirement that $m_{n} = o\{n^{1/(2\alpha+2)}\}$ comes from the following reason.
In the proof of Theorem \ref{thm: main1}, we require that there exists a sufficiently small constant $c>0$ such that, with probability approaching one,  $|\hat{\kappa}_{k} - \kappa_{\ell}| \geq c |\kappa_{k}-\kappa_{\ell}|$ for all $1 \leq k \leq m_{n}$ and $\ell \neq k$. 
Since $| \hat{\kappa}_{k} - \kappa_{\ell}| \geq |\kappa_{k} - \kappa_{\ell}| - | \hat{\kappa}_{k} - \kappa_{k}|$ and  $\sup_{k \geq 1}| \hat{\kappa}_{k}-\kappa_{k}| \leq ||| \hat{K} - K ||| = O_{\Prr}(n^{-1/2})$ by Lemma 4.2 in \cite{Bo00}, it suffices to have that $n^{1/2} \inf_{1 \leq k \leq m_{n}, \ell \neq k} | \kappa_{k} - \kappa_{\ell}| \to \infty$. Now, for any $1 \leq k \leq m_{n}$ and $\ell \neq k$, $| \kappa_{k} - \kappa_{\ell} | \geq \min \{ \kappa_{k} - \kappa_{k+1}, \kappa_{k-1} - \kappa_{k} \} \geq C_{1}^{-1}k^{-\alpha-1}\geq C_{1}^{-1} m_{n}^{-\alpha-1}$, and to ensure that $n^{1/2}m_{n}^{-\alpha-1} \to \infty$, we need that $m_{n} = o\{n^{1/(2\alpha+2)}\}$. In addition, in order that $m_{n} \sim n^{1/(\alpha+2\beta)}$ satisfies $m_{n}=o\{n^{1/(2\alpha+2)}\}$, we need that $\beta > \alpha/2+1$. 

The theorem shows that the value of $\gamma$ does not affect rates of convergence of $\hat{b}$, which is perhaps not surprising in view of the definition of $\hat{b}$. What is interesting is the fact that $\hat{b}$ with $m_{n}$ properly chosen is rate optimal, which means that smoothness of $b(s,t)$ in $s$ does not affect difficulty of estimating $b$. This is in sharp contrast with nonparametric estimation of a bivariate regression function. It should be noted that the results of Theorem \ref{thm: main1} continue to hold even if the condition that $|b_{j,k}| \leq C_{1}j^{-\gamma}k^{-\beta}$ for all $j,k \geq 1$ is replaced by a weaker condition that $| b_{j,k} | \leq C_{1} \ell_{j}k^{-\beta}$ for all $j,k \geq 1$ for some (given) positive sequence $\{ \ell_{j} \}_{j=1}^{\infty}$ such that $\sum_{j=1}^{\infty} \ell_{j}^{2} < \infty$. However, the value of $\gamma$ does matter for the analysis of the second estimator $\tilde{b}$. 


\cite{CrMa13} study prediction based on the estimator $\hat{b}$. They prove that, assuming $\Ep\{Y(t)\} = \Ep \{X(t)\} = 0$ for all $t \in I$, the estimator $\hat{Y}_{n+1} (s) = \int_{I} \hat{b}(s,t) X_{n+1}(t) dt$ with an appropriate choice of the cut-off level $m_{n}$ attains the minimax rate for estimation of $\Ep( Y_{n+1} \mid X_{n+1} )$ under the mean integrated squared error (MISE). 
Importantly, the prediction problem considered in \cite{CrMa13} is related to but substantially different from the problem of estimating $b$ considered in the present paper; the former is not an ill-posed inverse problem (is not a type of problems formulated as solving an integral equation; cf. Remark \ref{rem: ill-posed inverse problem}), and \cite{CrMa13} do not derive sharp rates of convergence for $\hat{b}$ itself and hence do not cover Theorem \ref{thm: main1}
 (the proof of Theorem 2 in \cite{CrMa13} does not lead to the results of our Theorem \ref{thm: main1} since from the beginning their proof is bounding $\Ep[\| \int_{I} \{ \hat{b}(\cdot,t) - b(\cdot,t) \}X_{n+1}(t) dt \|^{2}]$, and \cite{CrMa13} assume a stronger moment condition on $\xi_{k}$; see (6) in their paper). 
 \cite{PaQi12} and \cite{HoKi15} analyze the estimator $\hat{b}$ with dependent functional data, but they only prove consistency of $\hat{b}$ and thus do not cover Theorem \ref{thm: main1}. Precisely speaking, they prove consistency of the integral operator with kernel $\hat{b}$ for the operator norm.

Next, we derive rates of convergence for our second estimator. 

\begin{theorem}
\label{thm: main2}
Consider the estimator $\tilde{b}$ defined in (\ref{eq: second estimator}). Suppose that Assumption \ref{as: 1} is satisfied. Furthermore, suppose that $\gamma > \beta/2+1$. Then provided that $\max \{ m_{n,1}, m_{n,2} \} = o\{n^{1/(2\alpha+2)}\}$, we have 
\begin{equation}
||| \tilde{b} - b |||^{2} = O_{\Prr} \left \{ n^{-1} (m_{n,1} + m_{n,2}^{\alpha+1}) +m_{n,1}^{-2\gamma+1} + m_{n,2}^{-2\beta + 1} \right \}. \label{eq: bound1}
\end{equation}
Therefore, by choosing $m_{n,1} \sim \min \{ n^{1/(2\gamma)}, (n/\log n)^{1/(2\alpha+2)} \}$ and $m_{n,2} \sim n^{1/(\alpha + 2\beta)}$, 
we have 
\begin{align}
&||| \tilde{b} - b |||^{2}   = O_{\Prr} \left [ \max \{ (n/\log n)^{-(2\gamma -1)/(2\alpha +2)}, n^{-(2\gamma - 1)/(2\gamma)}, n^{-(2\beta-1)/(\alpha + 2\beta)} \} \right ].  \label{eq: bound2}
\end{align}
\end{theorem}

\begin{figure}[htbp]
 \begin{minipage}{0.24\hsize}
  \begin{center}
   \includegraphics[width=0.99\hsize]{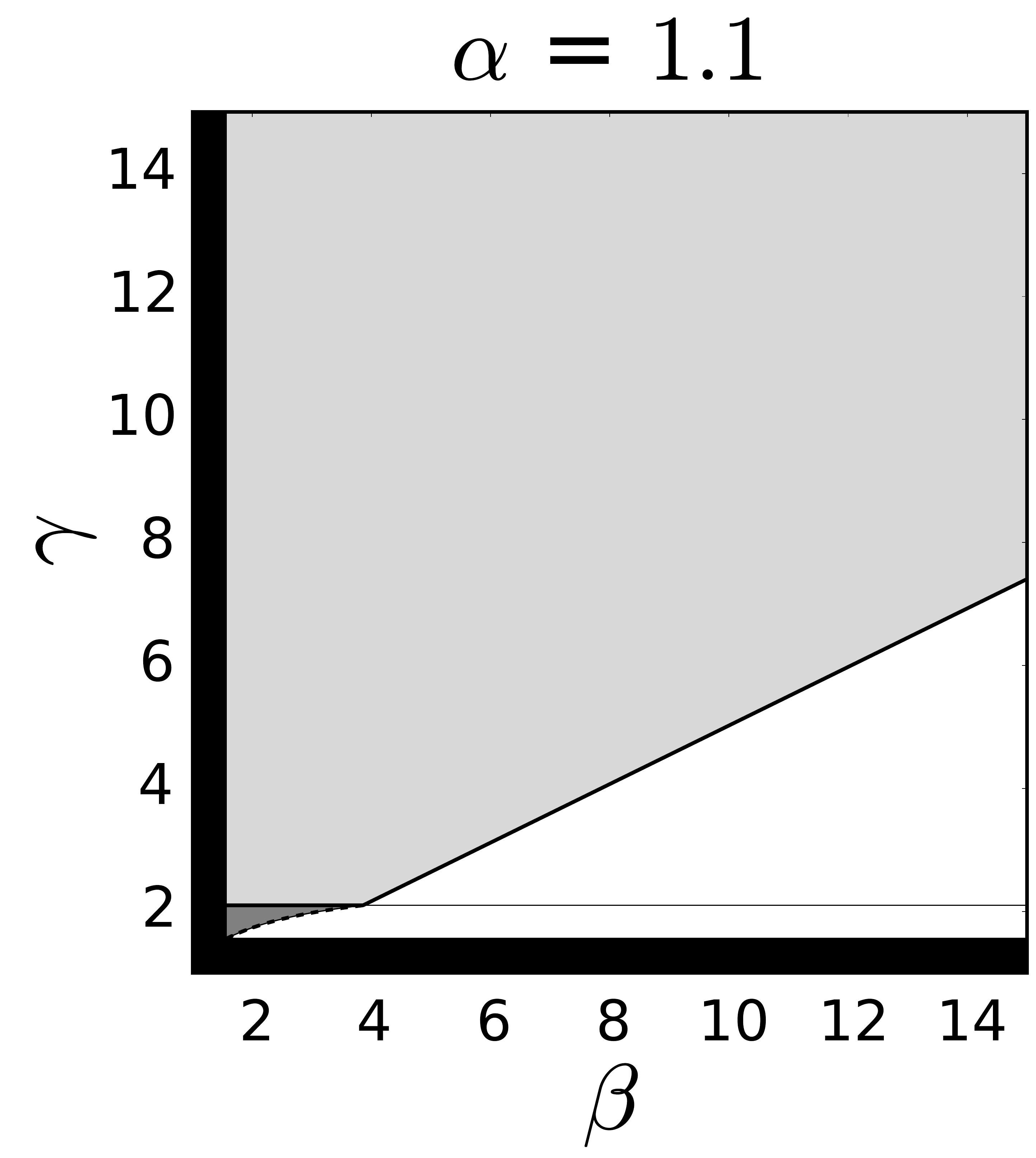}
  \end{center}
  \label{fig:one}
 \end{minipage}
 \begin{minipage}{0.24\hsize}
 \begin{center}
  \includegraphics[width=0.99\hsize]{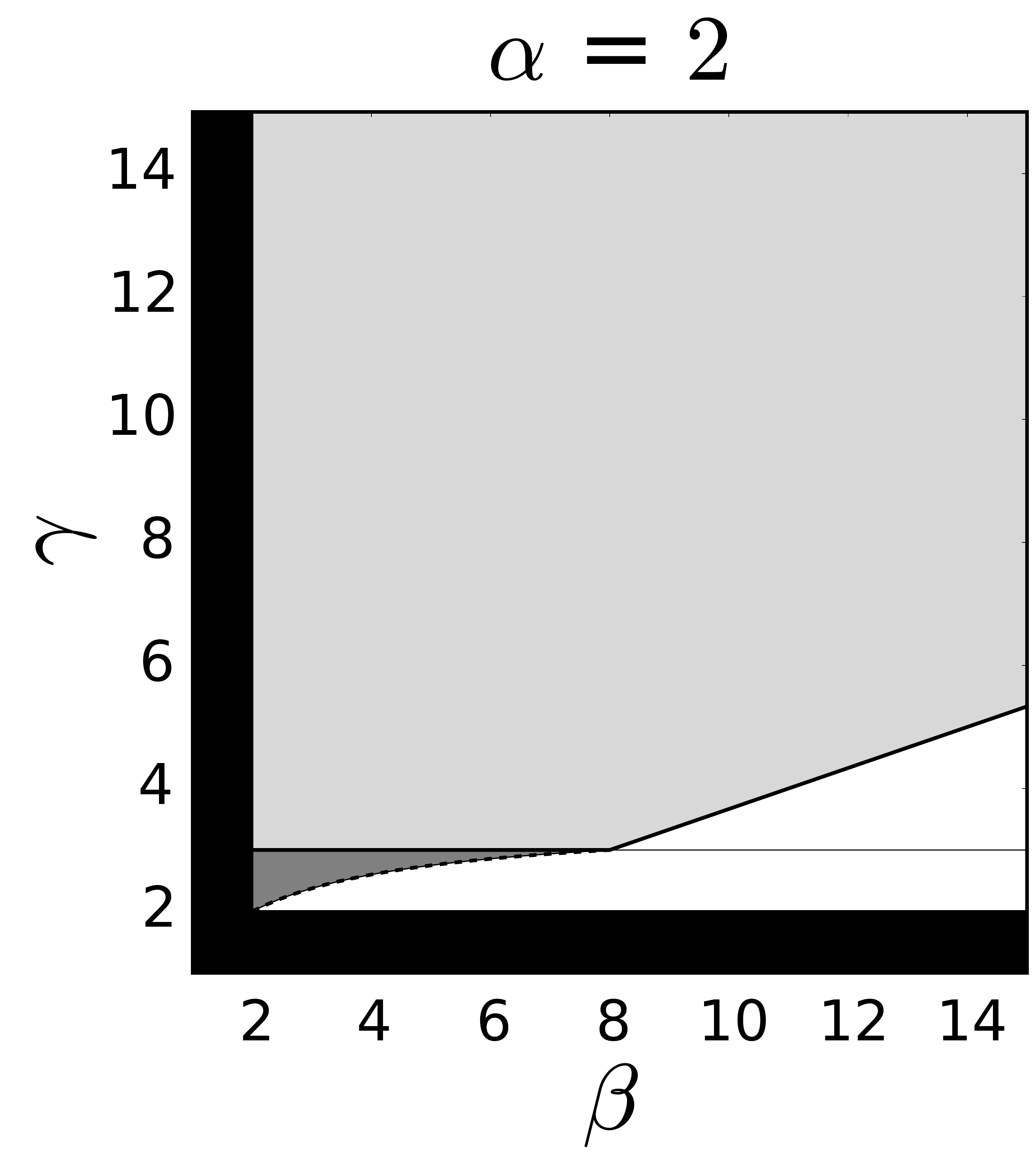}
 \end{center}
  \label{fig:two}
 \end{minipage}
 \begin{minipage}{0.24\hsize}
 \begin{center}
  \includegraphics[width=0.99\hsize]{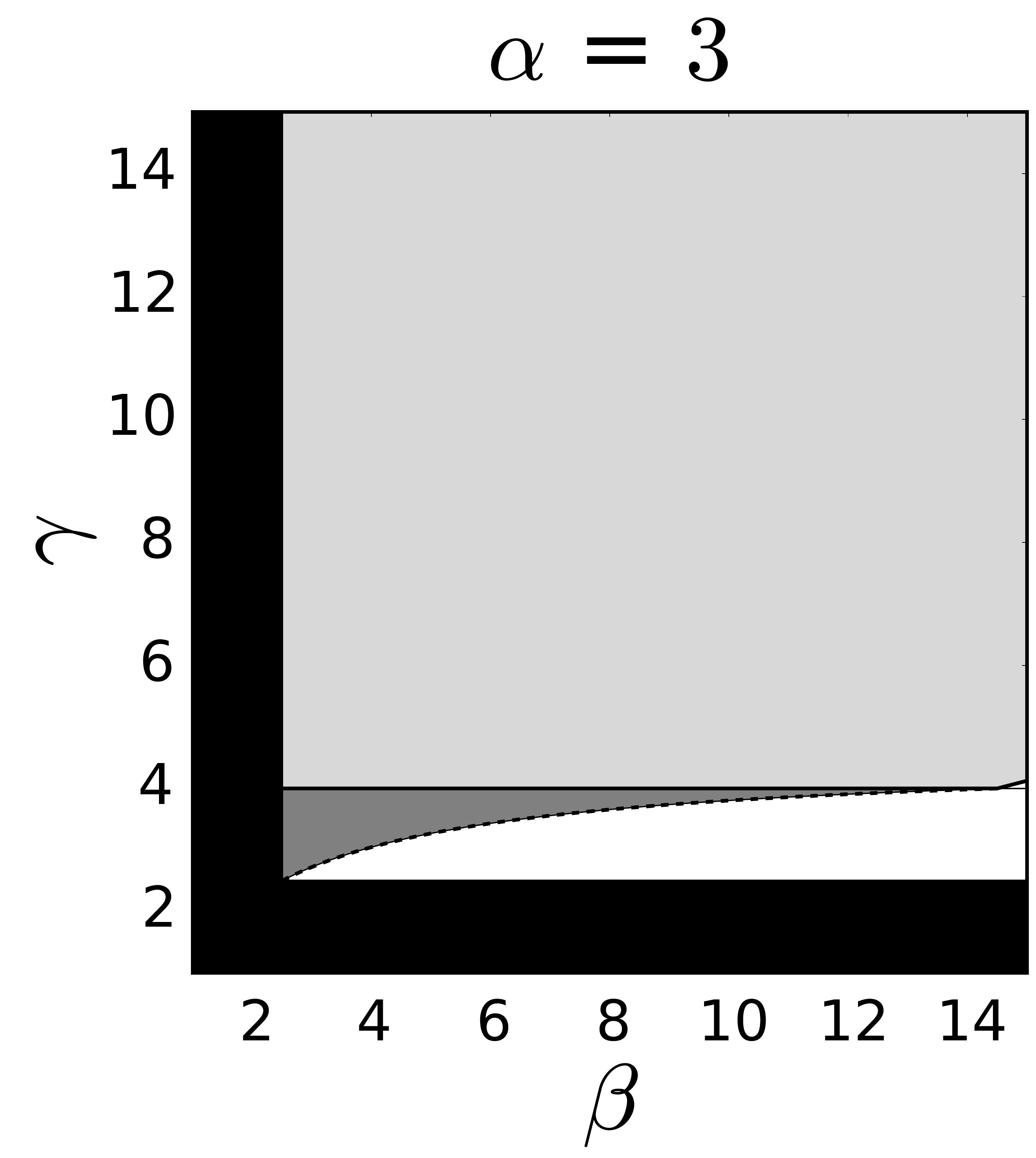}
 \end{center}
  \label{fig:three}
 \end{minipage}
 \begin{minipage}{0.24\hsize}
  \begin{center}
   \includegraphics[width=0.99\hsize]{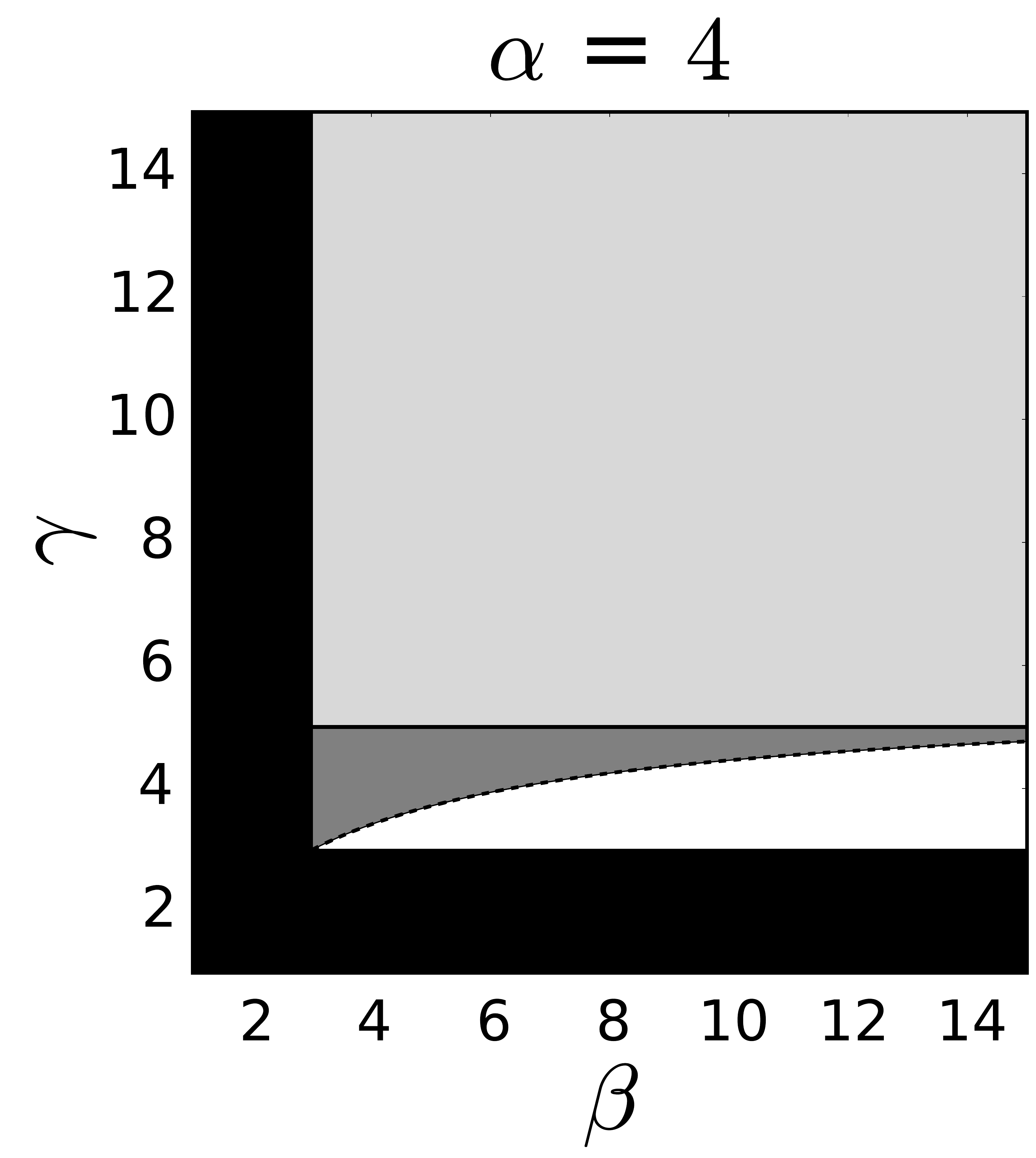}
  \end{center}
  \label{fig:four}
 \end{minipage}
 \caption{Regions of $(\beta,\gamma)$ for different values of $\alpha$.
When the parameters $(\beta,\gamma)$ are contained in the light gray region $(A)$ or the dark gray region $(B)$, the estimator $\tilde{b}$ attains the rate $n^{-(2\beta-1)/(\alpha+2\beta)}$. The black region corresponds to the region where $\beta \leq \alpha/2+1$ or $\gamma \leq \alpha/2+1$. \label{fig:parameters}}
\end{figure}

Since the estimator $\tilde{b}(s,t)$ depends on $\hat{\phi}_{1}(s),\dots,\hat{\phi}_{m_{n,1}}(s)$, accumulation of these estimation errors contributes to the term $n^{-1}m_{n,1}$ in the bound (\ref{eq: bound1}), while the term $m_{n,1}^{-2\gamma+1}$ comes from the bias. Because of these terms,  $\gamma$ appears in the bound (\ref{eq: bound2}), and in contrast to $\hat{b}$, the second estimator $\tilde{b}$ has suboptimal rates in some cases (of course there could be a room to improve upon the bound (\ref{eq: bound1})). 
Still, the estimator $\tilde{b}$ is able to attain the optimal rate $n^{-(2\beta - 1)/(\alpha + 2\beta)}$ provided that 
\begin{equation}
\beta 
\begin{cases}
\leq \frac{(2\gamma - 1)\alpha + 2\gamma}{2} & \text{if} \ \gamma > \alpha+1 \\
< \frac{(2\gamma -1)\alpha+2\alpha+2}{2(2\alpha - 2\gamma+3)} & \text{if} \ \gamma \leq \alpha+1
\end{cases}
,
\label{eq: restriction}
\end{equation}
which actually covers wide regions of $(\alpha,\beta,\gamma)$. 
Figure \ref{fig:parameters} depicts regions of $(\beta,\gamma)$ where $\tilde{b}$ attains the rate $n^{-(2\beta-1)/(\alpha + 2\beta)}$ for different values of $\alpha$.
We plot two regions $(A)=\{ (\beta,\gamma) :\alpha/2+1 < \beta \leq \{ (2\gamma - 1)\alpha + 2\gamma \}/2, \gamma > \alpha+1 \}$ and $(B) = \{ (\beta,\gamma) : \alpha/2+1 < \beta < \{ (2\gamma -1)\alpha+2\alpha+2 \}/\{2(2\alpha - 2\gamma+3)\}, \alpha/2+1 < \gamma \leq \alpha+1 \}$ in Figure \ref{fig:parameters}.

\cite{YaMuWa05} consider an estimator for $b$ that is related to but still different from our second estimator $\tilde{b}$. 
Their estimator is based on applying the functional PCA to both $X$ and $Y$. Let $L(s,t) = \Cov \{Y(s),Y(t)\}$ be the covariance function of $Y$, and let $L(s,t) = \sum_{j} \rho_{j} \psi_{j}(s) \psi_{j}(t)$ be the spectral expansion of $L$ where $\rho_{1} \geq \rho_{2} \geq \cdots > 0$ and $\{ \psi_{j} \}_{j=1}^{\infty}$ is an orthonormal basis of $L^{2}(I)$ (we assume here that this expansion is possible). Then $Y(s) = \Ep\{Y(s)\}+ \sum_{j} \zeta_{j}\psi_{j}(s)$ where $\zeta_{j}=\int_{I} [ Y(s)-\Ep\{Y(s)\} ] \psi_{j}(s)dt$, and observe that $b$ can be expanded in $L^{2}(I^{2})$ as $b(s,t) = \sum_{j,k} \{ \Cov (\zeta_{j},\xi_{k})/\kappa_{k} \} \psi_{j}(s) \phi_{k}(t)$.
The method of estimation of $b$ in \cite{YaMuWa05} is to approximate the infinite series $\sum_{j,k}$ by a finite series, and replace $\Cov (\zeta_{j},\xi_{k}), \kappa_{k},\psi_{j}$, and $\phi_{k}$ by their estimators. However, \cite{YaMuWa05} do not explicitly derive rates of convergence of this estimator, although it should be noted that \cite{YaMuWa05} assume that only discrete measurements with measurement errors for $X$ and $Y$ are available.
The analysis of the estimator of \cite{YaMuWa05} requires a substantially different set of assumptions than ours and thus is not pursued in the present paper. 

\subsection{Minimax lower bounds}
In this subsection, we derive minimax lower bounds for estimation of $b$. To this end, it is without loss of generality to narrow a class of distributions of $(X,Y)$, and we consider the following setting. 
Let $\alpha > 1, \beta > 1/2, \gamma > 1/2$, and $C_{1} > 1$ be given constants. 
Let $\mE$ be an $L^{2}(I)$-valued Gaussian random variable such that $\Ep( \langle f, \mE \rangle ) = 0$ and $\Ep ( \langle f,\mE \rangle^{2} ) > 0$ for all $f \in L^{2}(I)$ with $\| f \|=1$ (recall that an $L^{2}(I)$-valued random variable $Z$ is said to be Gaussian if $\langle f,Z \rangle$ is normally distributed for each $f \in L^{2}(I)$). 
Let $R(s,t) = \Ep \{ \mE (s) \mE (t) \}$ be the covariance function of $\mE$, and let $R(s,t) = \sum_{j=1}^{\infty} \lambda_{j} \phi_{j}(s) \phi_{j}(t)$ be the spectral expansion of $R$ where $\lambda_{1} \geq \lambda_{2} \geq \cdots > 0$ and $\{ \phi_{j} \}_{j=1}^{\infty}$ is an orthonormal basis of $L^{2}(I)$. Now, let $X = \sum_{k=1}^{\infty}k^{-\alpha/2} U_{k} \phi_{k}$ for $U_{1},U_{2},\dots$ being  independent uniform random variables on $[-3^{1/2},3^{1/2}]$ independent from $\mE$, and generate, as an $L^{2}(I)$-valued random variable, $Y(\cdot) = \int_{I} b(\cdot,t) X(t)dt + \mE (\cdot)$,
where $b \in L^{2}(I^{2})$. Since $U_{k}$  has mean zero and unit variance, we have $\kappa_{k}=k^{-\alpha}$, and so $
\kappa_{k} - \kappa_{k+1} = \alpha \int_{k}^{k+1} u^{-\alpha-1} du  \geq \alpha (k+1)^{-\alpha-1} \geq \alpha 2^{-\alpha-1} k^{-\alpha-1}$.
In addition, $\xi_{k}=k^{-\alpha/2} U_{k}$, and so $\Ep (\xi_{k}^{4}) = 9 k^{-2\alpha}/5$. Define 
\[
\mB (\beta,\gamma,C_{1}) = \left \{ b = \sum_{j,k} b_{j,k} \phi_{j} \otimes \phi_{k} : |b_{j,k}| \leq C_{1} j^{-\gamma} k^{-\beta}, \ \text{for all} \ j,k \geq 1 \right \}
\]
as a class of functions for $b$. 
\begin{theorem}
\label{thm: lower bound}
Work with the setting described as above. 
Then there exists a constant $c > 0$ such that $
\liminf_{n \to \infty} \inf_{\overline{b}^{n}} \sup_{b \in \mB (\beta,\gamma,C_{1})} 
\Prr_{b} \{   ||| \overline{b}^{n} - b |||^{2} \geq c n^{-(2\beta-1)/(\alpha + 2\beta)} \} > 0$, 
where $\Prr_{b}$ denotes the probability under $b$, and $\sup_{\overline{b}^{n}}$ is taken over all estimators $\overline{b}^{n}$ of $b$ based on $(X_{1},Y_{1}),\dots,(X_{n},Y_{n})$, independent copies of $(X,Y)$.
\end{theorem}

This theorem shows that, under Assumption \ref{as: 1}, the first PCA-based estimator $\hat{b}$ with $m_{n}$ properly chosen is minimax rate optimal, while the second PCA-based estimator $\tilde{b}$ with $(m_{n,1}, m_{n,2})$ properly chosen is minimax rate optimal provided that the additional restriction (\ref{eq: restriction}) is satisfied. 

\begin{remark}
One might be tempted to argue that the conclusion of Theorem \ref{thm: lower bound} would follow from the following observation: taking integration of $Y(s) = \int_{I} b(s,t) X(t)dt + \mathcal{E}(s)$, we arrive at the functional linear model 
\[
Y^{\sharp} = \int_{I} b_{\sharp}(t) X(t) + \varepsilon^{\sharp},
\]
where $Y^{\sharp} = \int_{I} Y(s) ds, b_{\sharp}(t) =\int_{I} b(s,t) ds$, and $\varepsilon^{\sharp} = \int_{I}\mathcal{E}(s) ds$. Since for any estimator $\overline{b}^{n}$ of $b$, $\| \overline{b}_{\sharp}^{n} - b_{\sharp} \| \leq ||| \overline{b}^{n} - b |||$ where $\overline{b}^{n}_{\sharp} (t)= \int_{I} \overline{b}^{n}(s,t) ds$, the conclusion of Theorem \ref{thm: lower bound} would follow from Theorem 1 in \cite{HaHo07}.
However, this argument contains a gap. The reason is that, when applying  Theorem 1 in \cite{HaHo07}, 
we implicitly restrict estimators of $b_{\sharp}$ to those based on $(Y^{\sharp}_{1},X_{1}),\dots,(Y_{n}^{\sharp},X_{n})$, thereby discarding the information that the entire paths of $Y_{1},\dots,Y_{n}$ are fully observed, which results in restricting a class of estimators of $b_{\sharp}$. Therefore, formally, the conclusion of Theorem \ref{thm: lower bound} does not directly follow from Theorem 1 in \cite{HaHo07}.
The proof of Theorem \ref{thm: lower bound} builds on constructing a suitable sequence of conditional distributions of $Y$ given $X$, and since $Y$ takes values in $L^{2}(I)$, we have to construct a sequence of distributions on $L^{2}(I)$, which is a significant difference from \cite{HaHo07}. To this end, we employ the theory of Gaussian measures on Banach spaces \cite[cf.][Chapter VIII]{St11}. 
\end{remark}

\section{Simulation results}
\label{sec: simulation}


In this section, we present simulation results to verify the performance of the estimators in the finite sample.
We consider the following data generating process. Let $\phi_{1} \equiv 1, \phi_{j+1}(t) = 2^{1/2} \cos (j\pi t) \ \text{for} \ j \geq 1$, and generate $(X,Y)$ as follows:
\begin{align*}
&Y(\cdot) = \int_{I} b(\cdot,t) X(t)dt + \mE (\cdot), \ X = \sum_{k=1}^{50} k^{-\alpha/2} U_{k} \phi_{k}, \ \mE = \sum_{j=1}^{50} j^{-1.1/2} Z_{j} \phi_{j}, \\
&b = \sum_{j,k=1}^{50} b_{j,k} \phi_{j} \otimes \phi_{k}, \ b_{1,1} = 0.3, \ b_{j,k} = 4(-1)^{j+k} j^{-\gamma} k^{-\beta} \ \text{for} \ (j,k) \neq (1,1), 
\end{align*}
where $U_{k} \sim \mathrm{Unif}. [-3^{1/2},3^{1/2}]$ and $Z_{j} \sim N(0,1)$ are all independent, and the following sample sizes for $n$ are examined: $400, 600, \dots, 2800, 3000$.
We consider the following configurations for $(\alpha, \beta, \gamma)$: 

\[
(1.2, 3, 2.5),  \ (1.2, 3, 3), \  (1.2, 3, 4), \ (2.4, 3, 2.5), \  (2.4, 3, 3), \  (2.4, 3, 4),
\] 
which verify the restriction (\ref{eq: restriction}). 
%
The number of repetitions for each simulation is $1000$. The numerical results obtained in this section were carried out by using the matrix language Ox \citep{Do02}. 
%

\begin{figure}[htbp]
 \begin{center}
 \includegraphics[width=0.99\hsize]{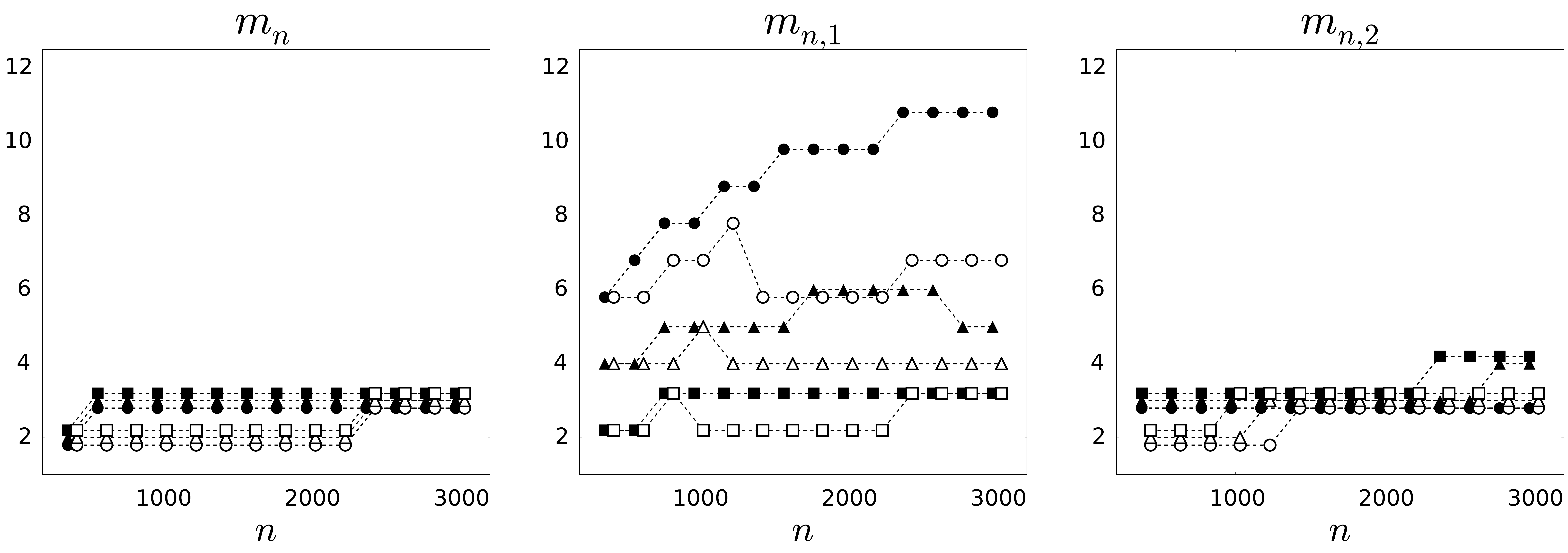}
 \end{center}
 \caption{The values of $m_{n}$ (left panel) and $(m_{n,1}, m_{n,2})$ (middle and right panels) minimizing MISE against $n$ for each parameter configuration. $(\alpha,\beta) = (1.2, 3.0)$ (black marker) and $(\alpha,\beta) = (2.4, 3.0)$ (white marker), and  $\gamma = 2.5$ (circle), $\gamma = 3.0$ (triangle) and $\gamma = 4.0$ (square).  \label{fig: mset}}
\end{figure}

\begin{figure}[htbp]
 \begin{center}
 \includegraphics[width=0.79\hsize]{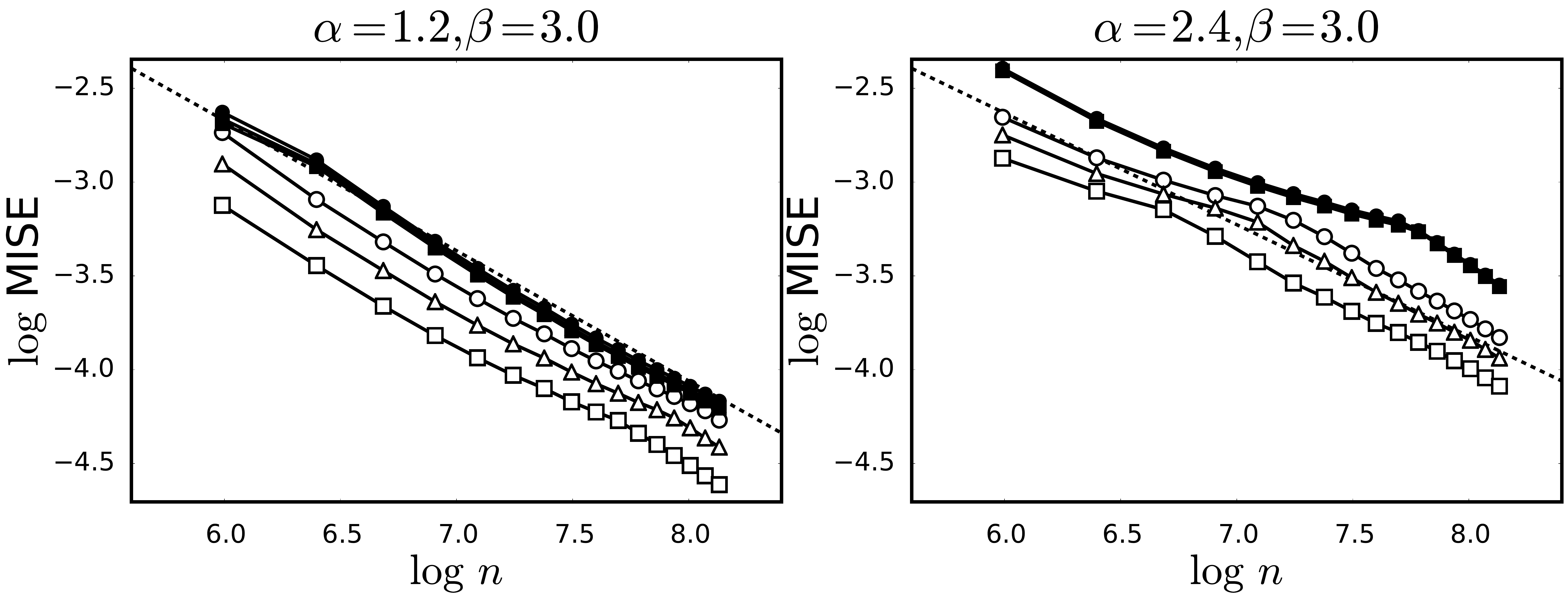}
 \end{center}
 \caption{Plots of the $\log$ MISE against $\log n$ for each parameter configuration : $(\alpha,\beta) = (1.2, 3.0)$ (left panel) and $(\alpha,\beta) = (2.4, 3.0)$ (right panel), and  $\gamma = 2.5$ (circle), $\gamma = 3.0$ (triangle) and $\gamma = 4.0$ (square). The black and white markers correspond to values of the $\log$ MISE of  $\hat{b}$ and $\tilde{b}$, respectively. The dashed line has slope $-(2\beta - 1) / (\alpha + 2 \beta)$. \label{fig: logn}}
\end{figure}

In this experiment, we simulate values of the MISE of $\hat{b}$ for $m_{n} \in \{1,\dots,20 \}$ and $\tilde{b}$ for $(m_{n,1},m_{n,2}) \in \{1,\dots,20 \}^{2}$ in each case, and report the optimal MISE. 
The selected values of $m_n$ and $(m_{n,1},m_{n,2})$ in each configuration are reported in Figure \ref{fig: mset}.
It is observed that 1) the values of $m_{n,1}$ selected become smaller as $\gamma$ increases, 2) the values of $m_{n,2}$ are less sensitive to $n$ than those of $m_{n,1}$, and 3) the values of $m_n$ are close to those of $m_{n,2}$.
Next, Figure \ref{fig: logn} plots the values of the log MISE against $\log n$. 
It is observed that 1) the values of the log MISE of $\hat{b}$ are almost identical for different values of $\gamma$; 
2) in contrast, the log MISE of $\tilde{b}$  decreases as $\gamma$ increases, but the slope is not sensitive to the value of $\gamma$, which indicates that the rate at which the MISE of $\tilde{b}$ decreases is independent of $\gamma$, but the constant depends on $\gamma$ and decreases as $\gamma$ increases;
3) all the slopes are close to $-(2\beta-1)/(\alpha+2\beta)$, at least for large $n$. These observations are consistent with our theoretical results.
Finally, in this limited experiment, the second estimator $\tilde{b}$ performs better than the first estimator $\hat{b}$, especially when $\gamma=4$; the difference in the $\log$ MISE is roughly $0.5$ in that case, which means that the MISE of $\hat{b}$ is $e^{0.5} \approx 1.65$ times that of $\tilde{b}$. 

\section{Real data analysis}
\label{sec: real data analysis}

\subsection{Working hours and income data}
\label{subsec: real data analysis}

We investigate the relation between the lifetime pattern of working hours and total income using data from  National Longitudinal Survey of Youth conducted by \cite{NLSY12}.
This is a major dataset in a field of human resources and consists of a sample from 12,686 American youth born between 1957 and 1964.
We use data of yearly working time (hour) and total net family income in a year from Round 1 (1979 survey year) to Round 25 (2012 survey year).

\begin{figure}[htbp]
\begin{minipage}{0.47\hsize}
 \begin{center}
 \includegraphics[width=0.99\hsize]{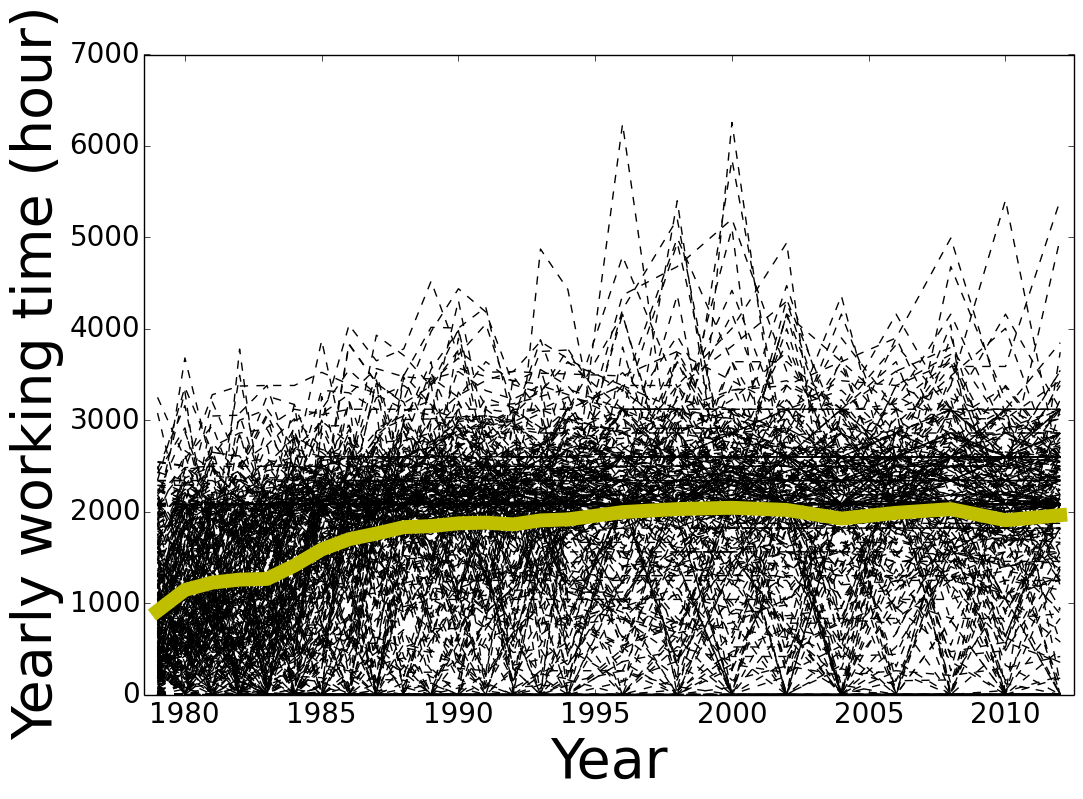}
 \end{center}
\end{minipage}
\begin{minipage}{0.47\hsize}
 \begin{center}
 \includegraphics[width=0.99\hsize]{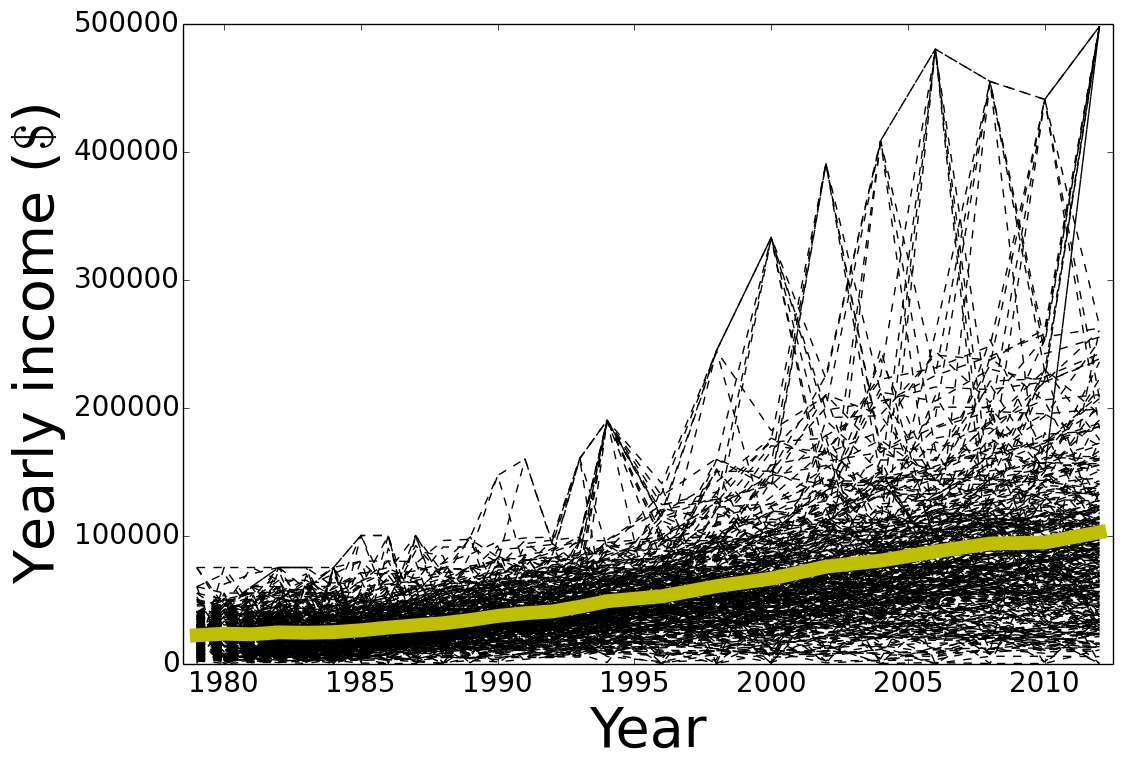}
 \end{center}
\end{minipage}
\caption{Plots of labor and income data. Yearly working hour (left panel) and yearly net income (right panel) against a survey year. Each black line is data for each cohort $i$ and yellow line is a mean function. \label{fig:data_labor}}
\end{figure}

We include cohorts who answer the question of all the 25 survey rounds and omit outliers who obtain more than $95\%$ quantiles of income.
Then we obtain working hour and income data of 353 observations and plot them in Figure \ref{fig:data_labor}.
In Figure \ref{fig:data_labor}, the black dashed lines show the working hour data $X_i(t)$ and the income data  $Y_i(t)$ for each respondent $i = 1,\ldots,353$.
The mean of working hour increases in the young age (about twenty to thirty) and the mean of income monotonically increases through all the rounds of the survey.
Since the income data is slightly discretized, some observations with high income take similar values.

We use the income as a response variable and the working hour as a predictor variable.
The values of $m_n$ and $(m_{n,1},m_{n,2})$ are selected by minimizing the cross-validation criteria as in \cite{YaMuWa05}: 
\begin{align*}
&\min_{m_{n}} \sum_{i=1}^n \int\left[  Y_i(s)  - \overline{Y}_{i}(s) - \int \hat{b}_{(-i,m_n)}(s,t)\{X_i(t) - \overline{X}(t)\}dt\right]^2ds, \\
&\min_{m_{n,1},m_{n,2}} \sum_{i=1}^n\int \left[  Y_i(s) - \overline{Y}_{i}(s) - \int \tilde{b}_{(-i,m_{n,1},m_{n,2})}(s,t)\{X_i(t) - \overline{X}(t)\}dt\right]^2ds,
\end{align*}
where $\hat{b}_{(-i,m)}(s,t)$ and $\tilde{b}_{(-i,m_1,m_2)}(s,t)$ are the estimates without $i$-th observation and with  the truncation levels $m_n$ and $(m_{n,1},m_{n,2})$, respectively. 
Using these criteria, we chose $m_n = 4$ for $\hat{b}$ and $(m_{n,1},m_{n,2}) = (4,4)$ for $\tilde{b}$.

\begin{figure}[htbp]
\begin{minipage}{0.47\hsize}
 \begin{center}
 \includegraphics[width=0.9\hsize]{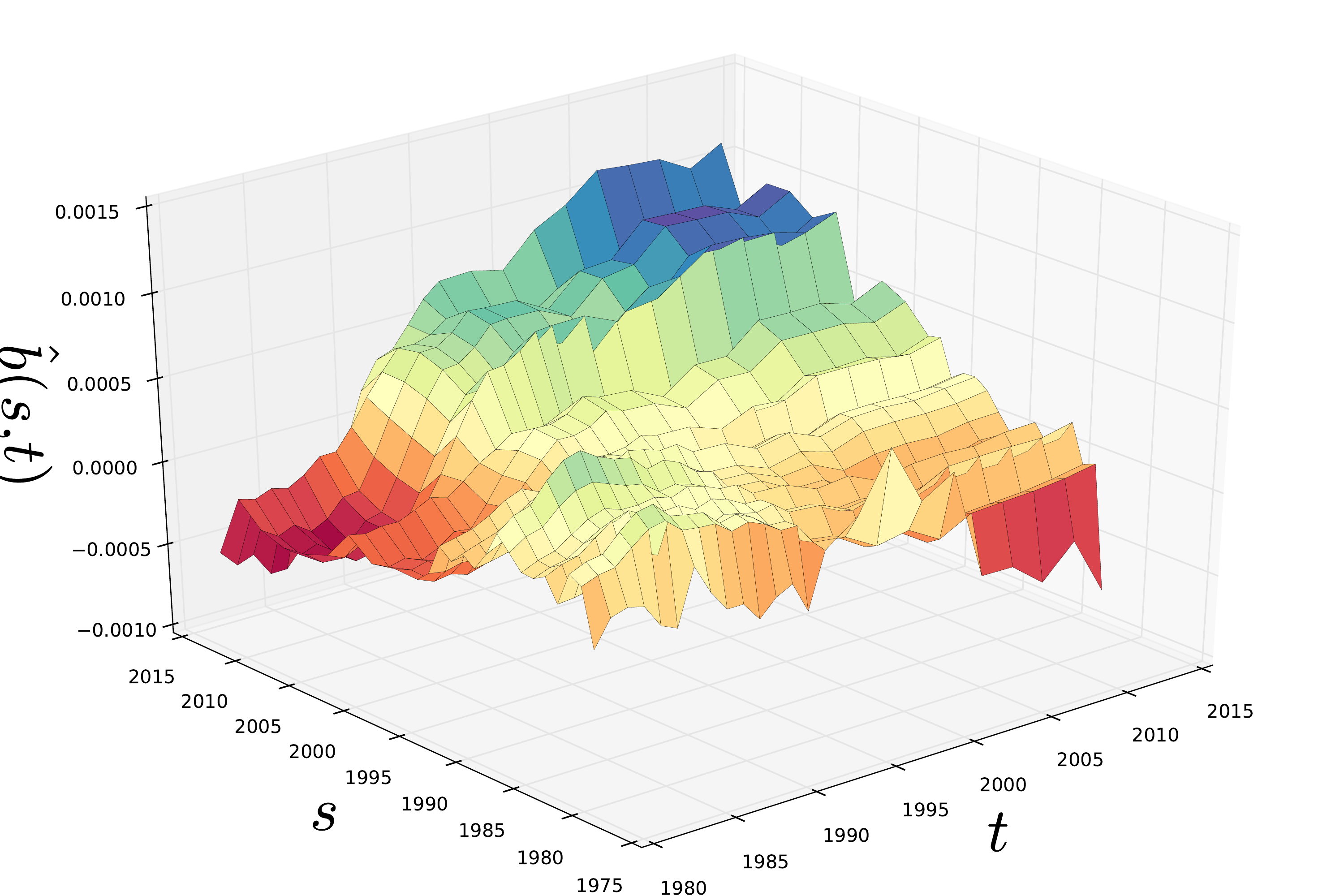} 
 \end{center}
 \caption{$\hat{b}(s,t)$ with the labor and income data. $m_n = 4$. \label{fig:labor_single}}
\end{minipage}
\begin{minipage}{0.47\hsize}
 \begin{center}
 \includegraphics[width=0.9\hsize]{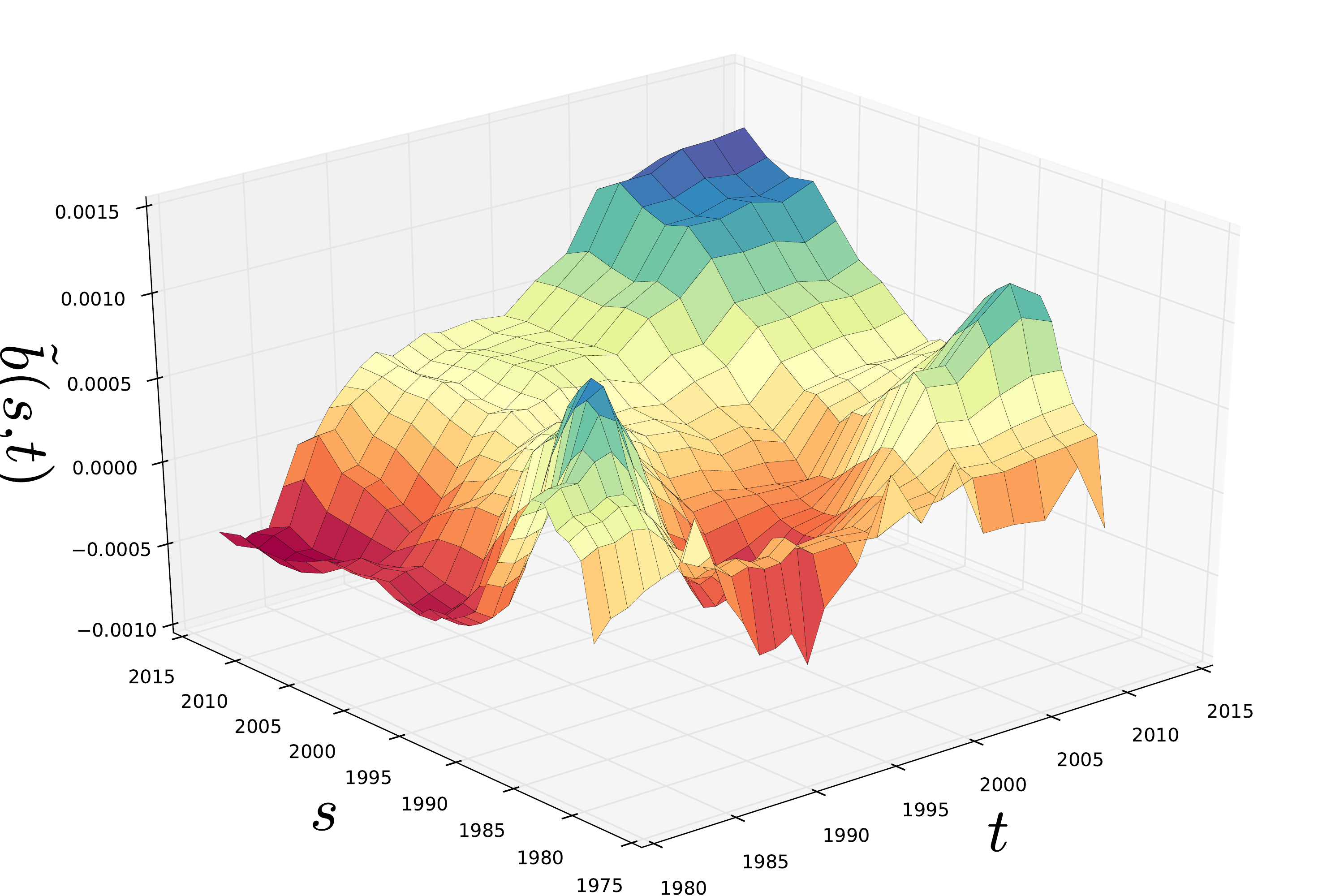} 
 \end{center}
 \caption{$\tilde{b}(s,t)$ with the labor and income data. $m_{n,1} =4$ and $ m_{n,2} = 4$.\label{fig:labor_double}}
\end{minipage}
\end{figure}

\begin{figure}[htbp]
\begin{minipage}{0.49\hsize}
\begin{minipage}{0.49\hsize}
 \begin{center}
 \includegraphics[width=0.99\hsize]{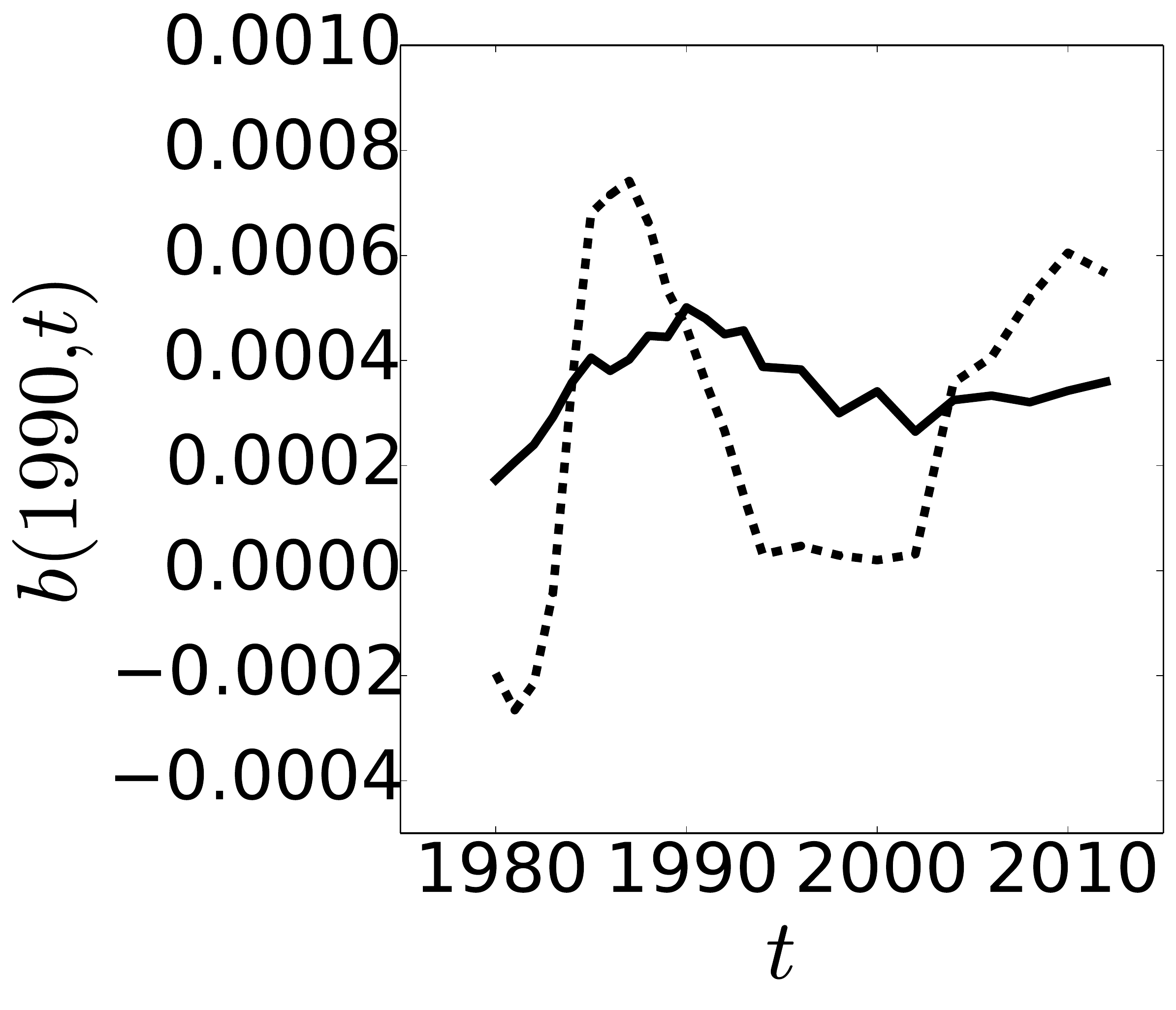}
 \end{center}
\end{minipage}
\begin{minipage}{0.49\hsize}
 \begin{center}
 \includegraphics[width=0.99\hsize]{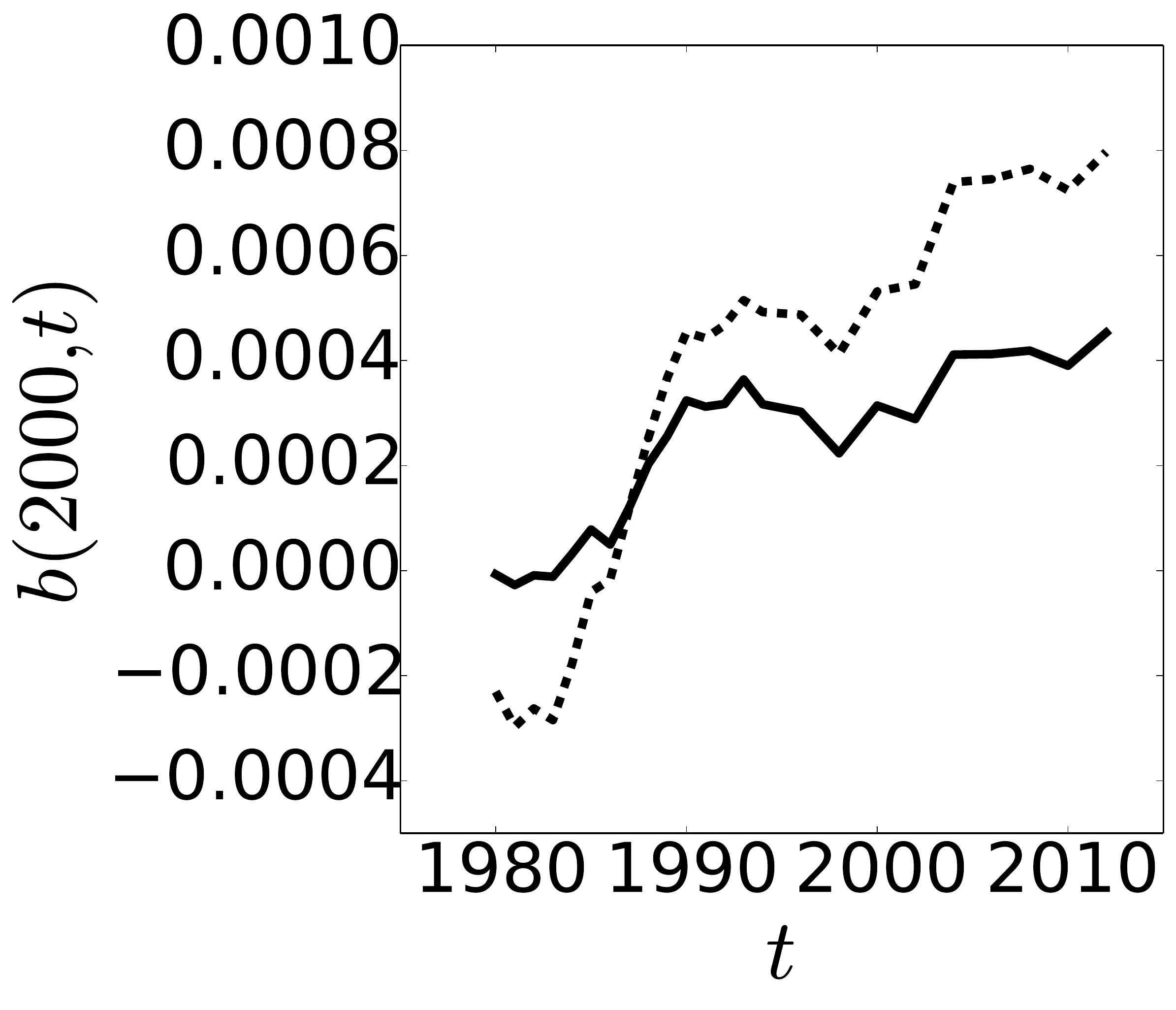}
 \end{center}
\end{minipage}
    \caption{Sliced $\hat{b}$ (solid) and $\tilde{b}$ (dashed) against $t$. \label{fig:labor_danmen1}}
\end{minipage}
\begin{minipage}{0.49\hsize}
\begin{minipage}{0.49\hsize}
 \begin{center}
 \includegraphics[width=0.99\hsize]{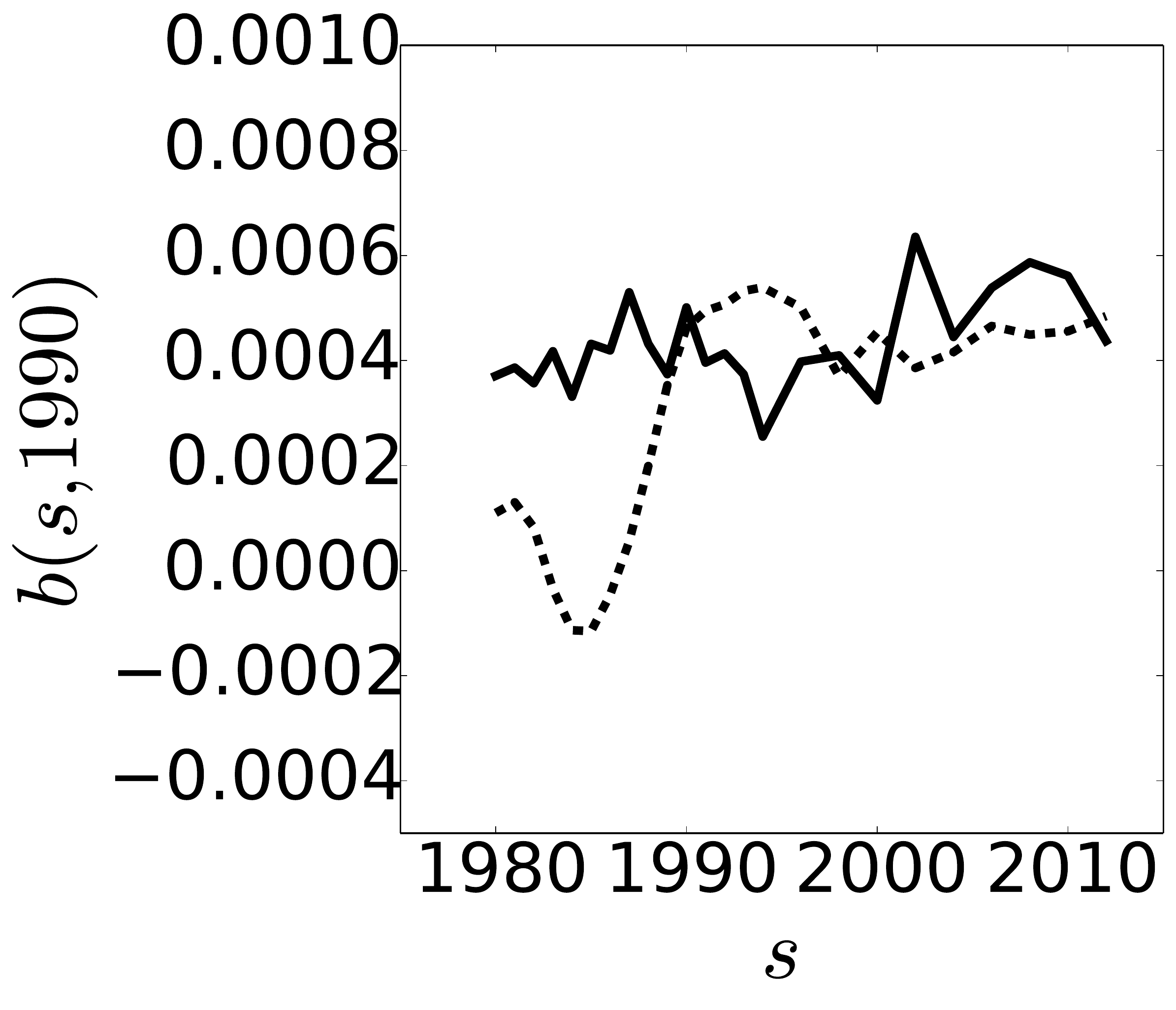}
 \end{center}
\end{minipage}
\begin{minipage}{0.49\hsize}
 \begin{center}
 \includegraphics[width=0.99\hsize]{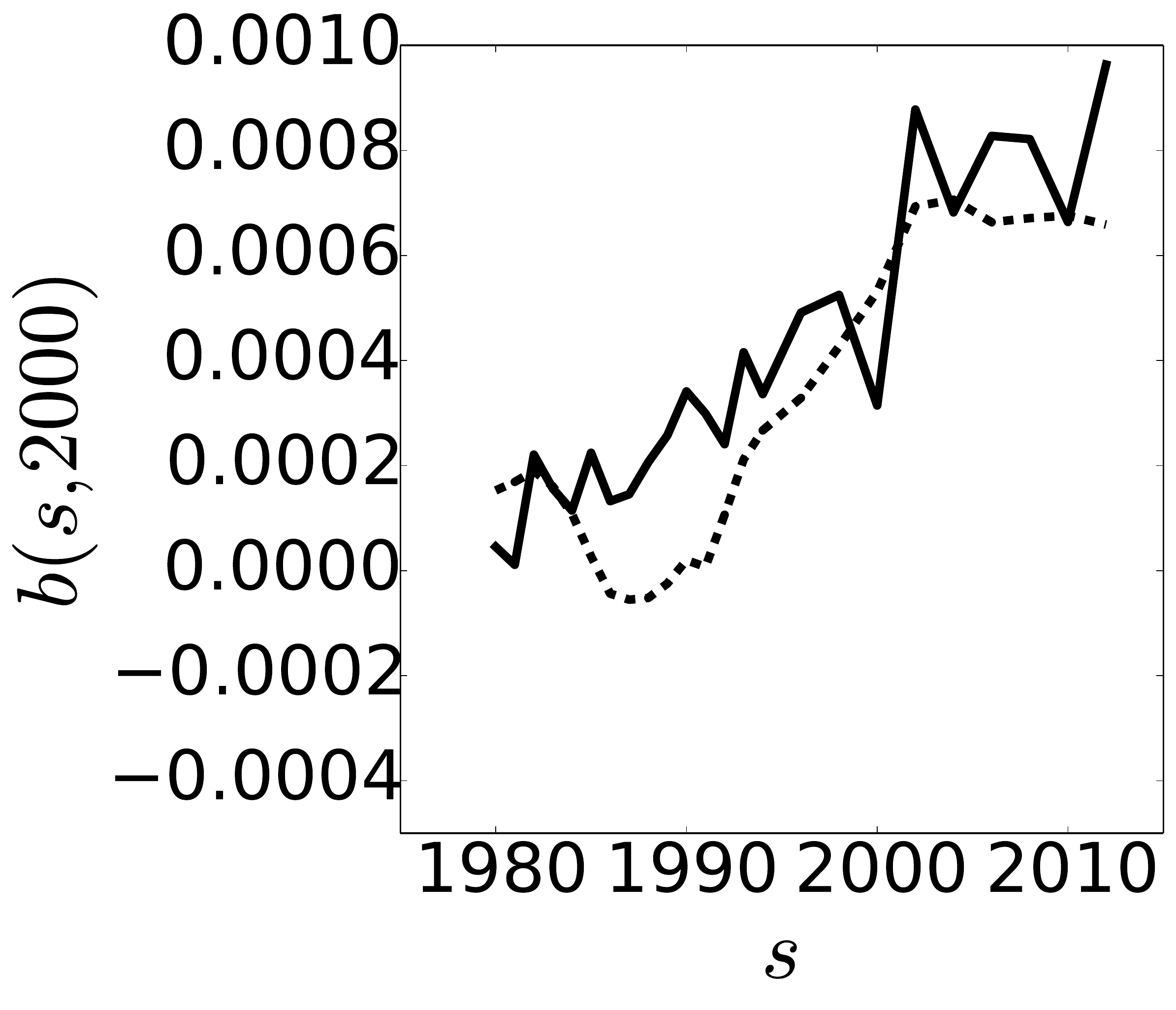}
 \end{center}
\end{minipage}
    \caption{Sliced $\hat{b}$ (solid) and $\tilde{b}$ (dashed) against $s$. \label{fig:labor_danmen2}}
\end{minipage}
\end{figure}

Figures \ref{fig:labor_single} and \ref{fig:labor_double} plot graphs of the estimates $\hat{b}$ and $\tilde{b}$, respectively. 
Figures \ref{fig:labor_danmen1} and \ref{fig:labor_danmen2} plot the slices of the estimates with $s= 1990, 2000$ and $t= 1990, 2000$. The overall shapes of the estimates as functions of $s$ or $t$ are roughly similar, but $\tilde{b}$ is smoother in $s$ than $\hat{b}$ because of the double truncation. 
Our functional regression analysis reveals that the working hour, not only in the advanced age but also in the middle age, can have positive effects on the income in the advanced age, and the positive effects get larger as the cohorts get older.
In contrast, the working hour in the young age has negative effects on the income in the middle and advanced ages.
This negative effect is interpreted as follows: cohorts who work much in their young age are not highly educated and they earn low income when they get older.

\subsection{Electricity prices} \label{subsec: real data analysis2}

We investigate the mechanism of electricity spot prices of the German power market traded at the European Energy Exchange (EEX).
In the German electricity market, the amount of renewable energy sources has a certain effect on the demand for the electricity because of the purchase guarantee, and the wind power infeed has the largest influence (a detailed discussion is found in \cite{Li13}).
With this background, we analyze how the wind power infeed affects the electricity price in the German power market.

\begin{figure}[htbp]
\begin{minipage}{0.47\hsize}
 \begin{center}
 \includegraphics[width=0.99\hsize]{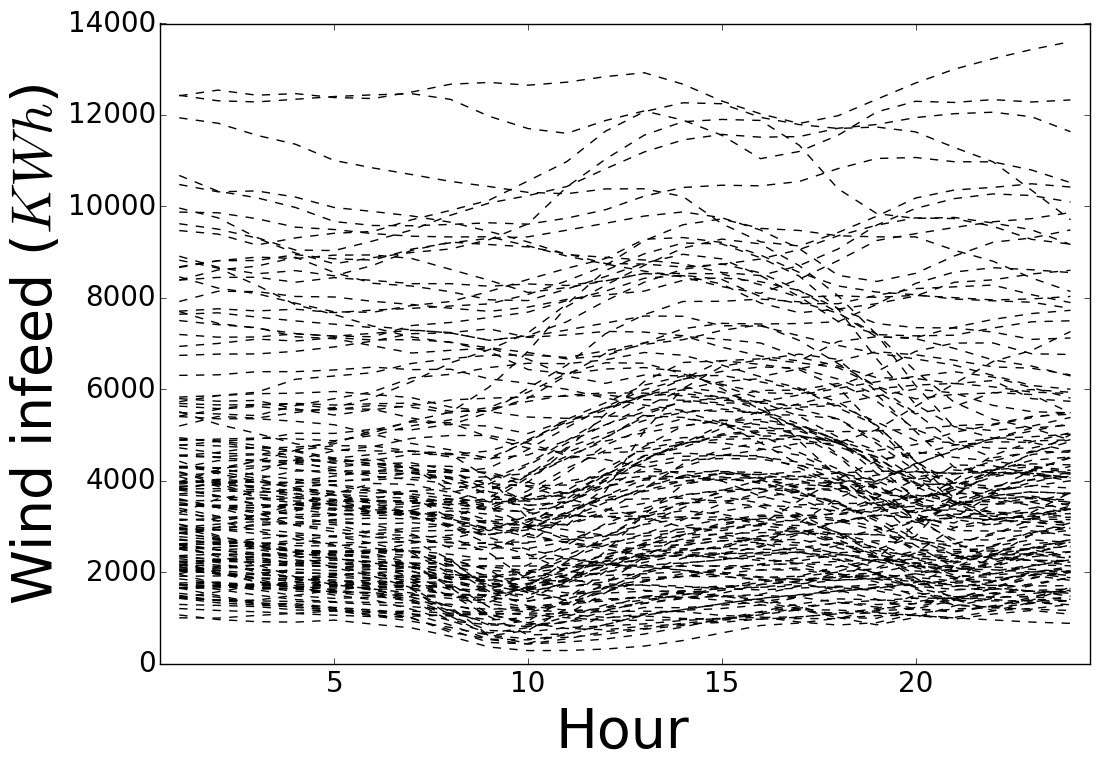}
 \end{center}
\end{minipage}
\begin{minipage}{0.47\hsize}
 \begin{center}
 \includegraphics[width=0.99\hsize]{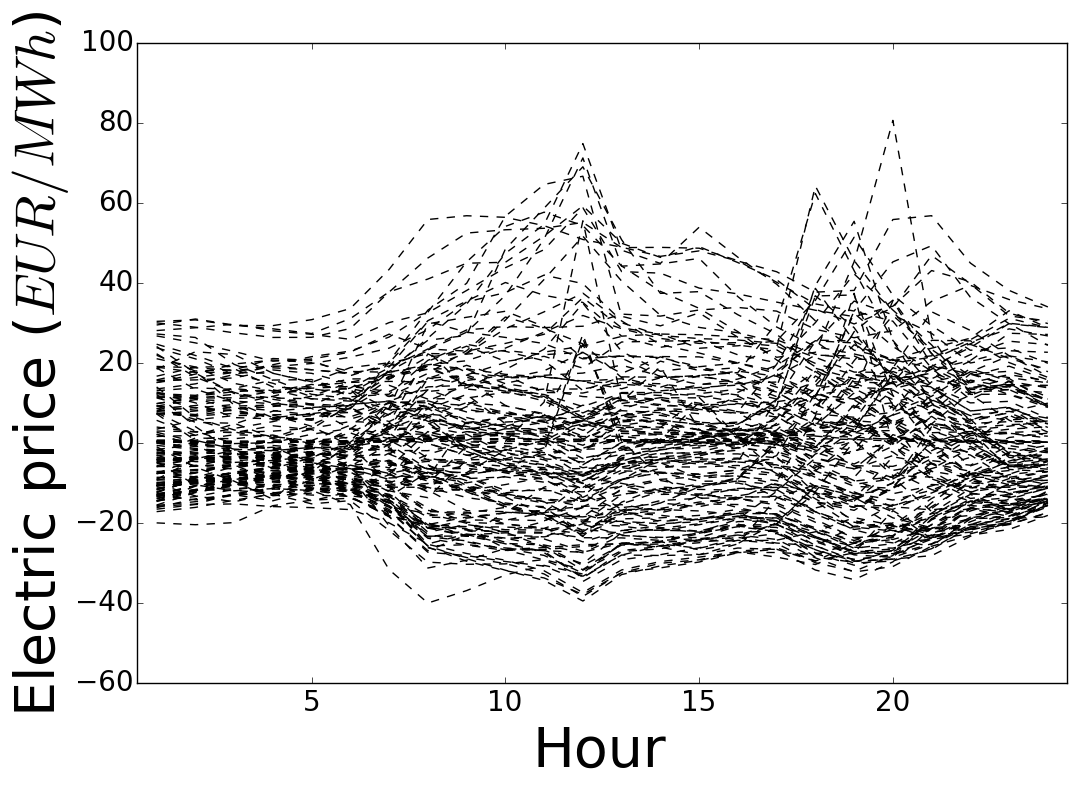}
 \end{center}
\end{minipage}
\caption{Plots of wind infeed and electricity price data. Wind power infeed (left panel) and centered electric price (right panel). \label{fig:data_electronic}}
\end{figure}

The data on prices of the German electricity market are taken from  European Energy Exchange, and the data on wind power in Germany are taken from the EEX Transparency Platform as in \cite{Li13}.
These data sets contain hourly electricity prices and wind power infeed from January 2006 to September 2008, and we take $Y_i(t)$ and $X_i(t)$ to be the electricity price and wind power infeed at time $t=1,\ldots,24$ and week $i =1,\ldots,143$  (each $Y_{i}(t)$ is centered around its sample mean).
 Figure \ref{fig:data_electronic} plots the data. Precisely speaking, the functional data in this example are likely to be dependent across $i$, but we expect that the convergence results in this paper could be extended to weakly dependent functional data. The formal analysis with dependent functional data is beyond the scope of the paper. 

\begin{figure}[H]
\begin{minipage}{0.47\hsize}
 \begin{center}
 \includegraphics[width=0.9\hsize]{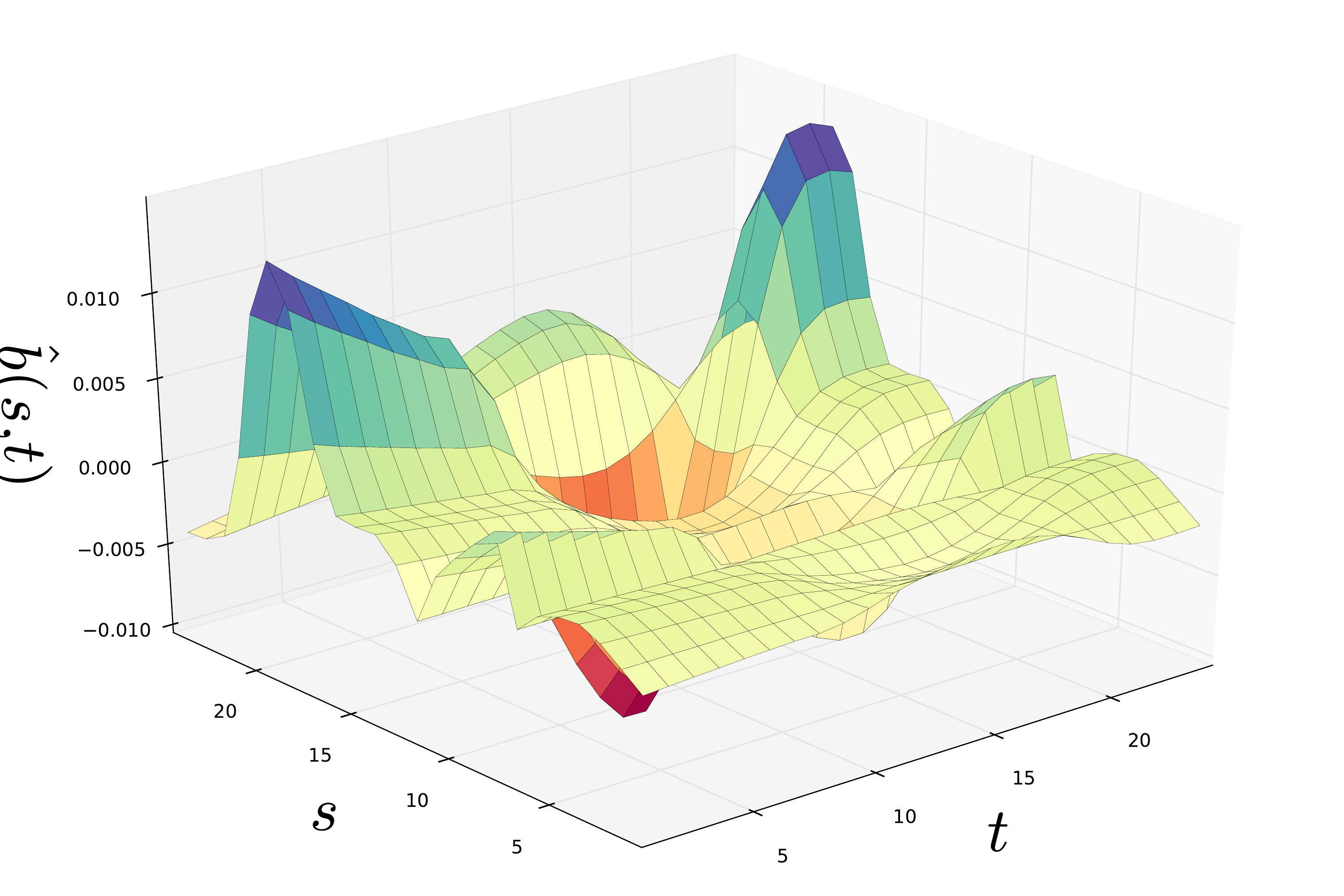}
 \end{center}
 \caption{$\hat{b}(s,t)$ with the electricity price and wind power infeed  data. $m_n = 2$. \label{fig:elec_single}}
\end{minipage}
\begin{minipage}{0.47\hsize}
 \begin{center}
 \includegraphics[width=0.9\hsize]{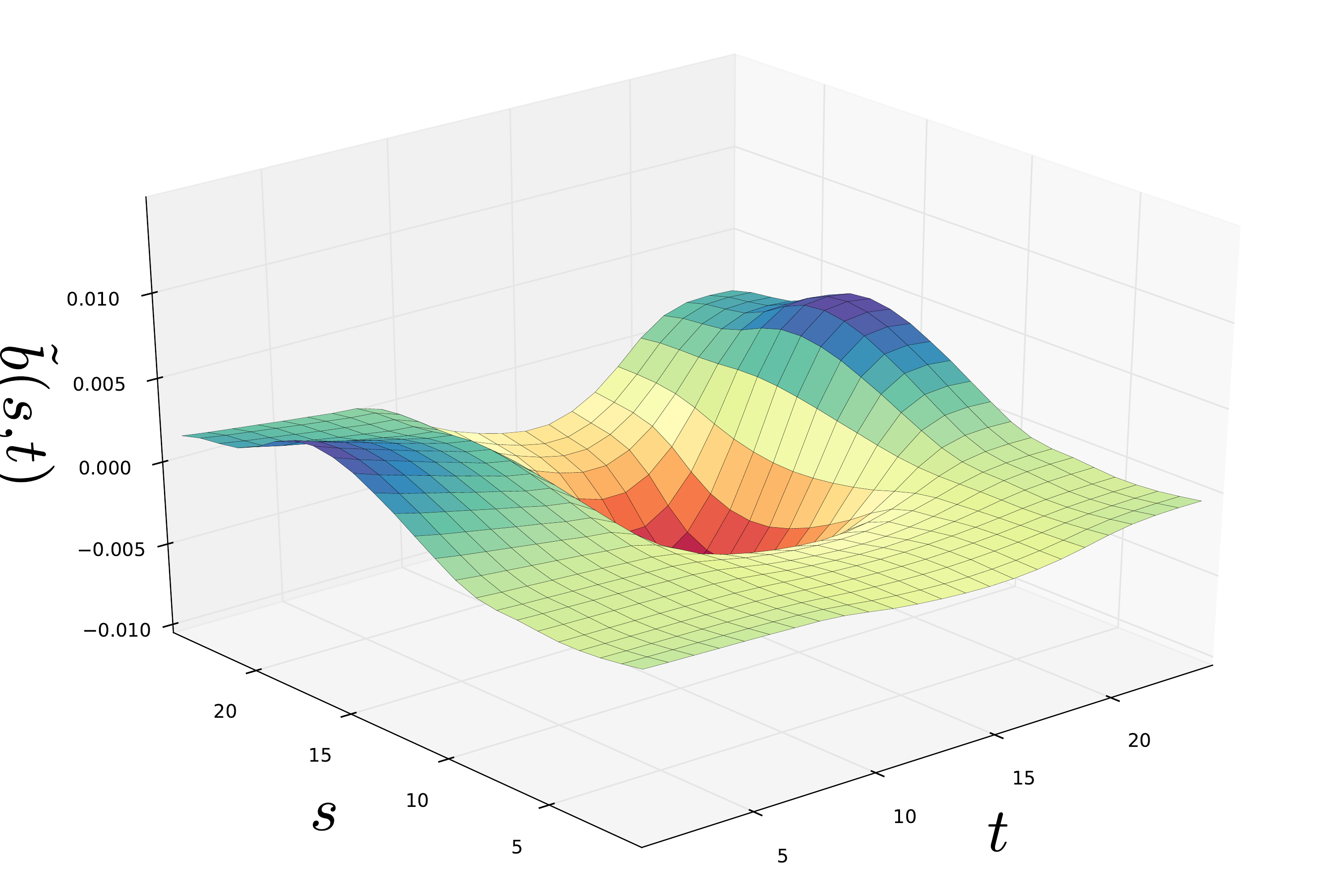}
 \end{center}
 \caption{$\tilde{b}(s,t)$ with the electricity price and wind power infeed data. $m_{n,1} =2$ and $ m_{n,2} = 1$.\label{fig:elec_double}}
\end{minipage}
\end{figure}

\begin{figure}[H]
\begin{minipage}{0.49\hsize}
\begin{minipage}{0.49\hsize}
 \begin{center}
 \includegraphics[width=0.99\hsize]{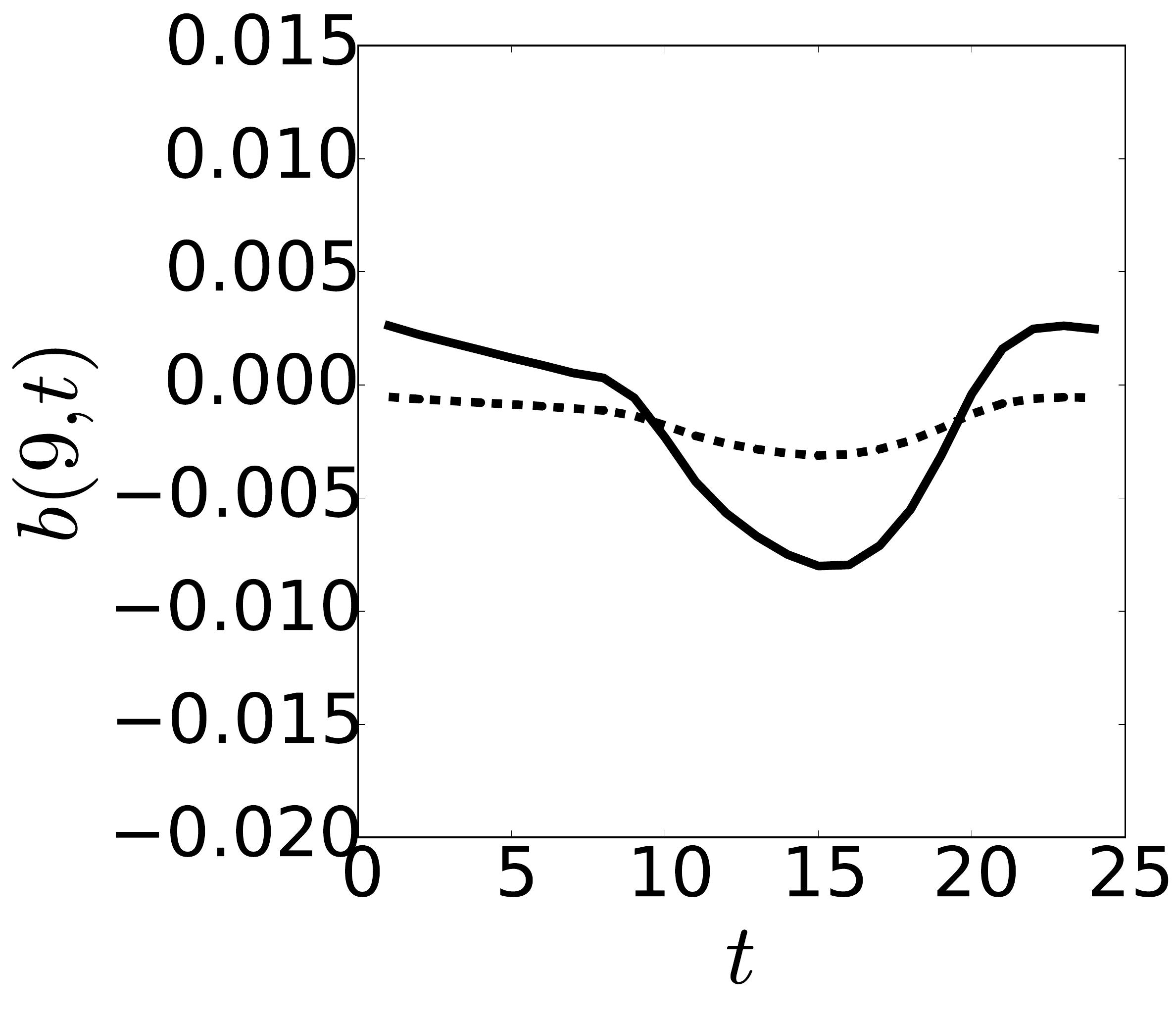}
 \end{center}
\end{minipage}
\begin{minipage}{0.49\hsize}
 \begin{center}
 \includegraphics[width=0.99\hsize]{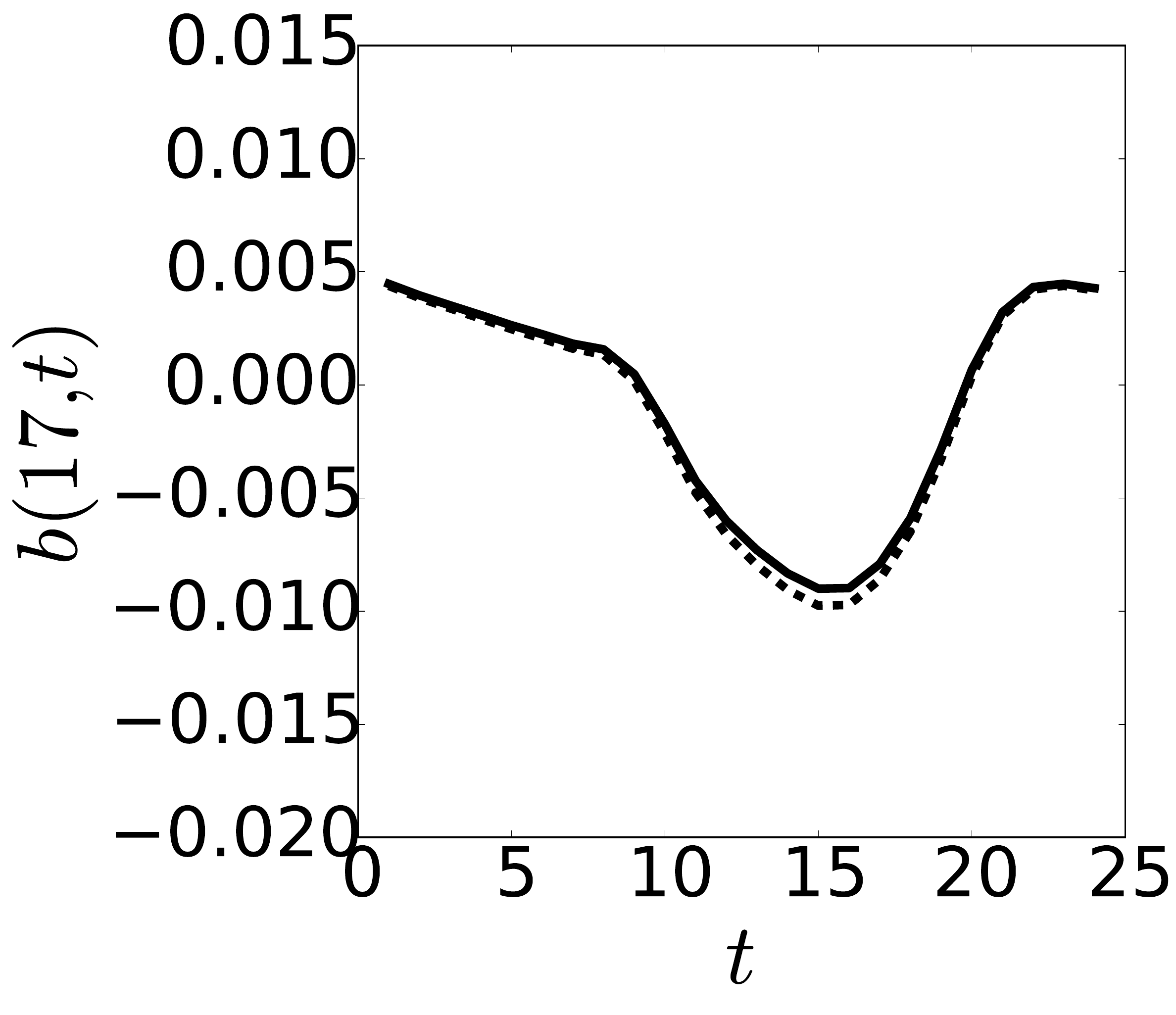}
 \end{center}
\end{minipage}
    \caption{Sliced $\hat{b}$ (solid) and $\tilde{b}$ (dashed) against $t$. \label{fig:elec_danmen1}}
\end{minipage}
\begin{minipage}{0.49\hsize}
\begin{minipage}{0.49\hsize}
 \begin{center}
 \includegraphics[width=0.99\hsize]{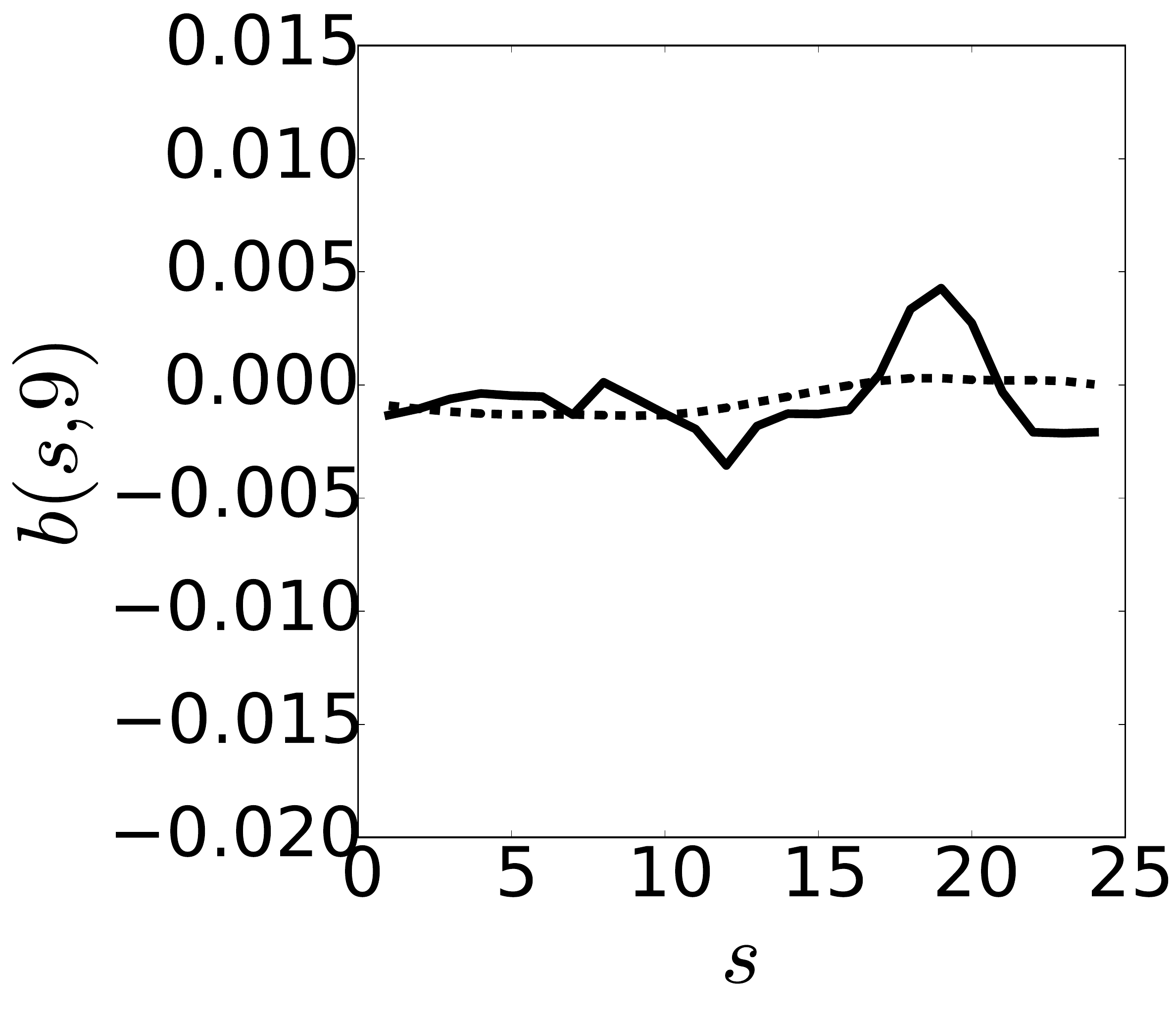}
 \end{center}
\end{minipage}
\begin{minipage}{0.49\hsize}
 \begin{center}
 \includegraphics[width=0.99\hsize]{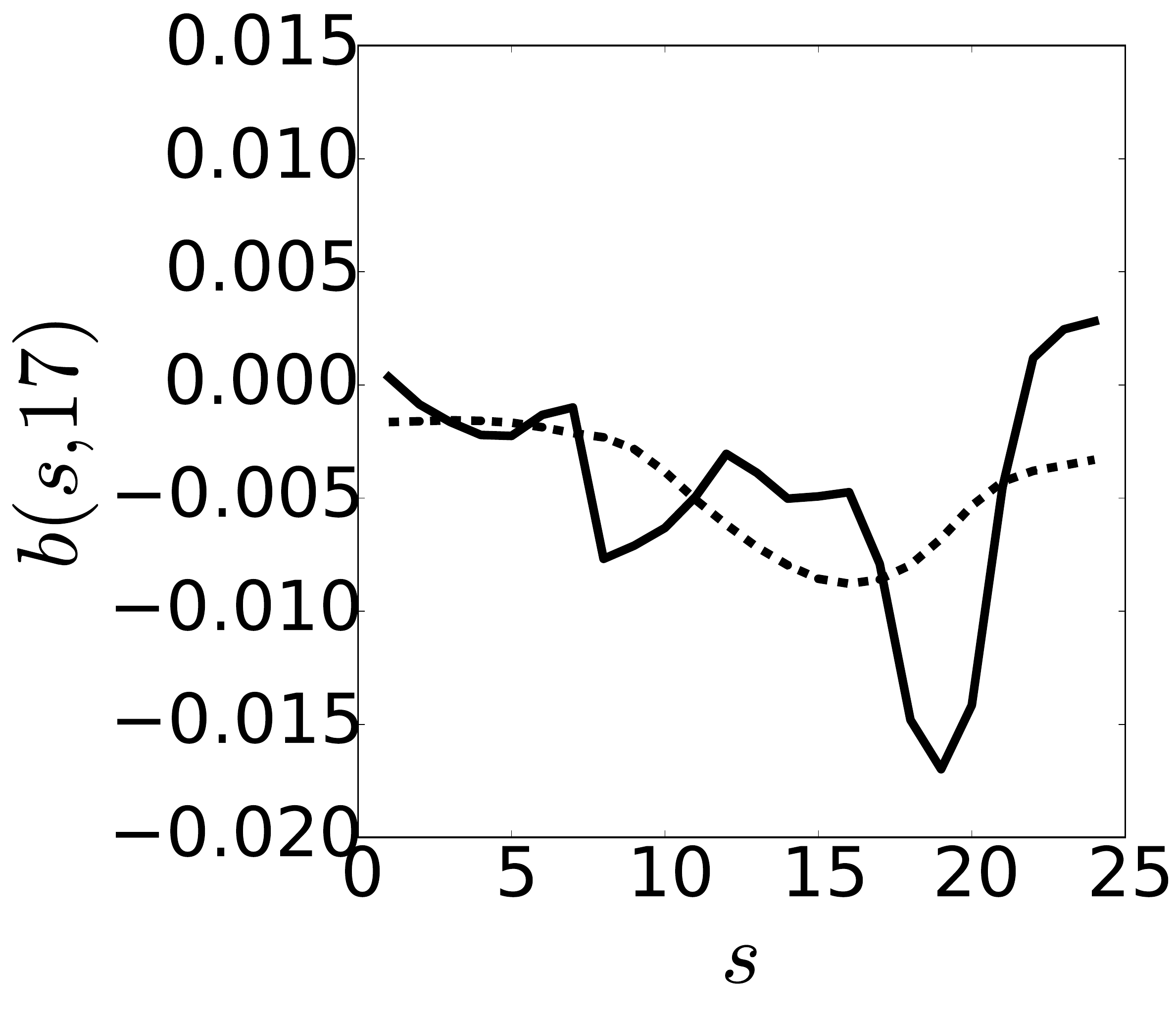}
 \end{center}
\end{minipage}
    \caption{Sliced $\hat{b}$ (solid) and $\tilde{b}$ (dashed) against $s$. \label{fig:elec_danmen2}}
\end{minipage}
\end{figure}

For this data set, we chose $m_n = 2$ for $\hat{b}$ and $(m_{n,1},m_{n,2}) = (2,1)$ for $\tilde{b}$ by the cross-validation. Figures \ref{fig:elec_single} and \ref{fig:elec_double} plot graphs of the estimates $\hat{b}$ and $\tilde{b}$, respectively.
Figures \ref{fig:elec_danmen1} and \ref{fig:elec_danmen2} plot slices of the estimates in the morning and evening (9 and 17 o'clock), which show that, as before, $\tilde{b}$ is smoother in $s$ than $\hat{b}$ because of the double truncation. 
Figure \ref{fig:elec_single} shows high fluctuations of the estimate $\hat{b}$, which make difficult to interpret the estimate. In contrast, from Figure \ref{fig:elec_double}, it is observed that $\tilde{b}$ is negative on the region $(s,t) \in [8,17] \times [8,17]$, but is close to zero on the other region. This shows that the wind power infeed has negative effects on the electricity price in the daytime, but except for the daytime, the effect of the wind power infeed is small.



\appendix

\section{Proof of Theorem \ref{thm: main1}}
\label{sec: proof}


In what follows, the notation $\lesssim$ signifies that the left hand side is bounded by the right hand side up to a constant that depends only on $\alpha,\beta,\gamma,C_{1}$. 
We first note that $\hat{b}$ is invariant with respect to choices of signs of $\hat{\phi}_{k}$'s, and so without loss of generality, we may assume that 
\begin{equation}
\int_{I} \hat{\phi}_{k}(t) \phi_{k}(t) dt \geq 0, \ \text{for all} \ k \geq 1. \label{eq: sign of phi}
\end{equation}
Recall that $m_{n} =  o\{ n^{1/(2\alpha+2)} \}$. Lemma 4.2 in \cite{Bo00} yields that 
\begin{equation}
\sup_{k \geq 1} | \hat{\kappa}_{k} - \kappa_{k} | \leq ||| \hat{K} - K |||.  \label{eq: eigenvalue bound}
\end{equation}
 Since 
\begin{align*}
\Ep (\| X - \Ep \{ X(\cdot) \} \|^{4}) &= \Ep \left \{ \left( \sum_{k=1}^{\infty} \xi_{k}^{2} \right )^{2} \right  \} = \sum_{k,\ell=1}^{\infty} \Ep(\xi_{k}^{2}\xi_{\ell}^{2})\\
& \leq \sum_{k,\ell=1}^{\infty}\{\Ep(\xi_{k}^{4})\}^{1/2} \{\Ep(\xi_{\ell}^{4})\}^{1/2} \lesssim \left ( \sum_{k=1}^{\infty} \kappa_{k} \right)^{2} \lesssim 1,
\end{align*}
we have that $||| \hat{K} - K ||| = O_{\Prr}(n^{-1/2})$. 
Define the event 
\[
A_{n} = \{ |\hat{\kappa}_{k}-\kappa_{\ell}| \geq  |\kappa_{k}-\kappa_{\ell}|/2 \ \text{for all} \ 1 \leq k \leq m_{n} \ \text{and for all}  \ \ell \neq k \}.
\]
It is seen that, since $
| \kappa_{k}-\kappa_{\ell} | \geq \min \{ \kappa_{k-1} - \kappa_{k}, \kappa_{k} - \kappa_{k+1} \} \geq C_{1}^{-1} k^{-\alpha-1} \geq C_{1}^{-1} m_{n}^{-\alpha-1}$
whenever $1 \leq k \leq m_{n}$ and $\ell \neq k$, and since $n^{-1/2} = o(m_{n}^{-\alpha-1})$, we have $
\Prr (A_{n}) \to 1$.
Furthermore, arguing as in \citet[][p.83-84]{HaHo07}, we have
\begin{equation}
\| \hat{\phi}_{k} - \phi_{k} \|^{2} \leq 8 \{1+o_{\Prr}(1) \} \hat{u}_{k}^{2}, \ \Ep( \hat{u}_{k}^{2} ) \lesssim k^{2}/n \label{eq: eigenfunction bound}
\end{equation}
where $o_{\Prr}(1)$ is uniform in $1 \leq k \leq m_{n}$. In what follows, we will freely use the estimates in (\ref{eq: eigenvalue bound}) and (\ref{eq: eigenfunction bound}). 
In particular, since $\beta  > 3/2$, we have 
\[
\sum_{k=1}^{m_{n}} k^{-2\beta} \| \hat{\phi}_{k} - \phi_{k} \|^{2} = O_{\Prr}\left (n^{-1} \sum_{j=1}^{\infty}k^{-2\beta + 2} \right)= O_{\Prr}(n^{-1}).
\]
In what follows, integrations such as $\int_{I} f(t)dt$ and $\iint_{I^{2}} R(s,t) ds dt$ are abbreviated as $\int f$ and $\iint R$. 

Let $\mE_{i} = Y_{i} - \Ep( Y_{i} \mid X_{i} )$, and expand $Y_{i}$ and $\mE_{i}$ as $Y_{i}  = \sum_{j} \eta_{i,j} \phi_{j}$ and $\mE_{i} = \sum_{j} \varepsilon_{i,j} \phi_{j}$,
where $\eta_{i,j}= \int Y_{i} \phi_{j}$ and $\varepsilon_{i,j} = \int \mE_{i} \phi_{j}$. 
Observe that $\hat{b}$ admits the following alternative expansion in $L^{2}(I)$: 
$\hat{b} = \sum_{j=1}^{\infty} \sum_{k=1}^{m_{n}} \hat{b}_{j,k} (\phi_{j} \otimes \hat{\phi}_{k})$,
where $\hat{b}_{j,k} = n^{-1} \sum_{i=1}^{n} \eta_{i,j}\hat{\xi}_{i,k}/\hat{\kappa}_{k}$. 
Since $\{ \phi_{j} \otimes \hat{\phi}_{k} \}_{j,k=1}^{\infty}$ is an orthonormal basis of $L^{2}(I^{2})$, expand $b$ as $b = \sum_{j,k} \breve{b}_{j,k} (\phi_{j} \otimes \hat{\phi}_{k})$ with $\breve{b}_{j,k} = \iint b (\phi_{j} \otimes \hat{\phi}_{k})$.
Now, setting $\eta_{i,j}^{c} = \eta_{i,j} - n^{-1}\sum_{i'=1}^{n}\eta_{i',j}$ and $\varepsilon_{i,j}^{c} =\varepsilon_{i,j} - n^{-1} \sum_{i'=1}^{n} \varepsilon_{i',j}$, observe that 
\[
\eta_{i,j}^{c} = \sum_{\ell} \breve{b}_{j,\ell} \hat{\xi}_{i,\ell} + \varepsilon_{i,j}^{c}.
\]
Plugging this expression into $\hat{b}_{j,k}$ together with the facts that 
\begin{equation}
n^{-1}\sum_{i=1}^{n} \hat{\xi}_{i,k} = 0, \ n^{-1}\sum_{i=1}^{n} \hat{\xi}_{i,\ell} \hat{\xi}_{i,k} = \iint \hat{K} (\hat{\phi}_{\ell} \otimes \hat{\phi}_{k}) 
= \begin{cases}
\hat{\kappa}_{k} & \text{if} \ \ell = k \\
0 & \text{if} \ \ell \neq k
\end{cases}
, \label{eq: orthogonality}
\end{equation}
we have $\hat{b}_{j,k}=\breve{b}_{j,k} + n^{-1} \sum_{i=1}^{n} \varepsilon_{i,j}\hat{\xi}_{i,k}/\hat{\kappa}_{k}$.
Therefore, 
\[
(\hat{b}_{j,k} - b_{j,k})^{2} \lesssim (\breve{b}_{j,k} - b_{j,k})^{2} + \hat{\kappa}_{k}^{-2} \left ( \frac{1}{n} \sum_{i=1}^{n} \varepsilon_{i,j}\hat{\xi}_{i,k} \right)^{2}.
\]
We divide the rest of the proof into three steps.

\textbf{Step 1}. We wish to bound $\sum_{j=1}^{\infty} \sum_{k=1}^{m_{n}} (\breve{b}_{j,k} - b_{j,k})^{2}$. We will use the following expansion: if $\inf_{\ell: \ell \neq k} | \hat{\kappa}_{k} - \kappa_{\ell} | > 0$, then
\begin{equation}
\hat{\phi}_{k} - \phi_{k} = \sum_{\ell: \ell \neq k} (\hat{\kappa}_{k}-\kappa_{\ell})^{-1} \phi_{\ell} \iint (\hat{K} - K) (\hat{\phi}_{k} \otimes \phi_{\ell}) + \phi_{k} \int (\hat{\phi}_{k}-\phi_{k}) \phi_{k}.
\label{eq: difference}
\end{equation}
See Lemma 5.1 in \cite{HaHo07}. Observe that
\begin{align*}
&\breve{b}_{j,k} - b_{j,k} 
=\sum_{\ell: \ell \neq k} b_{j,\ell} (\hat{\kappa}_{k}-\kappa_{\ell})^{-1} \iint (\hat{K}-K) (\hat{\phi}_{k} \otimes \phi_{\ell}) + b_{j,k} \int (\hat{\phi}_{k}-\phi_{k}) \phi_{k} \\
&=\sum_{\ell: \ell \neq k} b_{j,\ell} (\kappa_{k} - \kappa_{\ell})^{-1} \iint (\hat{K}-K) (\phi_{k} \otimes \phi_{\ell}) \\
&\quad + \sum_{\ell: \ell \neq k} b_{j,\ell} \{ (\hat{\kappa}_{k} - \kappa_{\ell})^{-1} - (\kappa_{k} - \kappa_{\ell})^{-1} \} \iint (\hat{K}-K) (\phi_{k} \otimes \phi_{\ell}) \\
&\quad + \sum_{\ell: \ell \neq k} b_{j,\ell} (\hat{\kappa}_{k} - \kappa_{\ell})^{-1} \iint (\hat{K}-K) \{ (\hat{\phi}_{k} - \phi_{k})\otimes \phi_{\ell} \}  + b_{j,k} \int (\hat{\phi}_{k}-\phi_{k}) \phi_{k} \\
&=: T_{j,k,1} + T_{j,k,2} + T_{j,k,3} + T_{j,k,4}.
\end{align*}

It is seen that $|T_{j,k,4}| \lesssim j^{-\gamma}k^{-\beta} \| \hat{\phi}_{k}-\phi_{k} \|$.
Next, since
\[
\left| \iint (\hat{K}-K) \{ (\hat{\phi}_{k} - \phi_{k})\otimes \phi_{\ell} \} \right| \leq ||| \hat{K}-K ||| \cdot \| \hat{\phi}_{k} - \phi_{k} \|,
\]
we have on the event $A_{n}$, 
\[
|T_{j,k,3}| \lesssim j^{-\gamma}||| \hat{K}-K ||| \cdot \| \hat{\phi}_{k} - \phi_{k} \| \sum_{\ell: \ell \neq k} \frac{\ell^{-\beta}}{|\kappa_{k}-\kappa_{\ell}|}.
\] 
In view of the assumption that $k^{-\alpha} \lesssim \kappa_{k} \lesssim k^{-\alpha}$, choose $k_{0} \geq 1$ and $C > 1$ large enough so that $\kappa_{k}/\kappa_{[k/C]} \leq 1/2$  and $\kappa_{[Ck]+1}/\kappa_{k} \leq 1/2$ for all $k \geq k_{0}$, where $[a]$ denotes the largest integer not exceeding $a$.  We may choose $k_{0}$ and $C$ in such a way that they depend only on $\alpha$ and $C_{1}$. Now, partition the sum $\sum_{\ell: \ell \neq k}$ into $\sum_{\ell=1}^{[k/C]}, \ \sum_{\ell=[k/C]+1,\neq k}^{[Ck]}$ and $\sum_{\ell =[Ck]+1}^{\infty}$.
Observe that 
\begin{align*}
\sum_{\ell=1}^{[k/C]} \frac{\ell^{-\beta}}{(\kappa_{\ell}-\kappa_{k})} &\leq \sum_{\ell=1}^{[k/C]} \frac{\ell^{-\beta}}{\kappa_{\ell}(1-\kappa_{k}/\kappa_{[k/C]})} \lesssim \sum_{\ell=1}^{[k/C]} \ell^{-\beta+\alpha} \\
&\lesssim 
\begin{cases}
1 & \text{if} \ \beta > \alpha + 1 \\
\log k & \text{if} \ \beta = \alpha + 1\\
k^{\alpha-\beta+1} & \text{if} \ \beta < \alpha + 1
\end{cases}
,
\end{align*}
\[
\sum_{\ell =[Ck]+1}^{\infty}\frac{\ell^{-\beta}}{(\kappa_{k}-\kappa_{\ell})} \leq \sum_{\ell = [Ck]+1}^{\infty} \frac{\ell^{-\beta}}{\kappa_{k} (1-\kappa_{[Ck]+1}/\kappa_{k})} \lesssim k^{\alpha} \sum_{\ell=[Ck]+1}^{\infty} \ell^{-\beta} \lesssim k^{\alpha-\beta+1}. 
\]
For $[k/C] < \ell < k$, observe that 
\begin{align*}
\kappa_{\ell} - \kappa_{k} &\geq k^{-\alpha-1}  \left \{ C_{1}^{-1} + k^{\alpha+1} \sum_{p=\ell}^{k-2} (\kappa_{p}-\kappa_{p+1}) \right \} \\
&\geq k^{-\alpha-1} C_{1}^{-1}\left \{ 1+\sum_{p=\ell}^{k-2} (k/p)^{\alpha+1} \right \} \\
&\geq k^{-\alpha-1}C_{1}^{-1} (k-\ell).
\end{align*}
Likewise, for $k < \ell \leq [Ck]$, $\kappa_{k} - \kappa_{\ell} \gtrsim k^{-\alpha-1} |k-\ell|$. Hence 
\[
\sum_{\ell=[k/C]+1,\neq k}^{[Ck]} \frac{\ell^{-\beta}}{|\kappa_{\ell}-\kappa_{k}|} \lesssim k^{\alpha+1} \sum_{\ell=[k/C]+1,\neq k}^{[Ck]} \frac{\ell^{-\beta}}{|k-\ell|} \lesssim k^{\alpha-\beta+1} \log k. 
\]
This yields that
\begin{equation}
\sum_{\ell: \ell \neq k} \frac{\ell^{-\beta}}{|\kappa_{k}-\kappa_{\ell}|} \lesssim 
\begin{cases}
1  & \text{if} \ \beta > \alpha + 1 \\
k^{\alpha-\beta+1} \log k & \text{if} \ \beta \leq \alpha+1
\end{cases}
,
\label{eq: delicate bound}
\end{equation}
and so on the event $A_{n}$, 
\[
|T_{j,k,3}| \lesssim j^{-\gamma} ( 1 + k^{\alpha-\beta+1} \log k ) \cdot ||| \hat{K}-K ||| \cdot \| \hat{\phi}_{k} - \phi_{k} \|.
\]

Turning to $T_{j,k,1}$, observe that for each $\ell \neq k$, 
\[
\iint (\hat{K}-K) (\phi_{k} \otimes \phi_{\ell}) = \frac{1}{n} \sum_{i=1}^{n} \xi_{i,k}\xi_{i,\ell} - \overline{\xi}_{k} \overline{\xi}_{\ell},
\]
which yields that $\Ep( T_{j,k,1}^{2} )$ is 
\begin{equation}
\lesssim n^{-1} \Ep \left \{ \xi_{k}^{2} \left ( \sum_{\ell: \ell \neq k} \frac{b_{j,\ell}}{\kappa_{k}-\kappa_{\ell}} \xi_{\ell} \right )^{2} \right \}   +\Ep \left \{ \overline{\xi}_{k}^{2} \left ( \sum_{\ell: \ell \neq k}  \frac{b_{j,\ell}}{\kappa_{k}-\kappa_{\ell}}\overline{\xi}_{\ell} \right )^{2} \right \}. \label{eq: T1}
\end{equation}
The first term on the right hand side is bounded by 
\[
n^{-1}\{\Ep(\xi_{k}^{4})\}^{1/2} \left [ \Ep \left \{ \left ( \sum_{\ell: \ell \neq k} \frac{b_{j,\ell} \xi_{\ell} }{\kappa_{k}-\kappa_{\ell}}\right )^{4} \right \} \right ]^{1/2},
\]
where $\{\Ep(\xi_{k}^{4})\}^{1/2} \lesssim k^{-\alpha}$. Now, observe that 
\begin{align*}
&\Ep \left \{ \left ( \sum_{\ell: \ell \neq k} \frac{b_{j,k} \xi_{\ell} }{\kappa_{k}-\kappa_{\ell}}\right )^{4} \right \} 
\leq 
\sum_{\ell_{1}: \ell_{1} \neq k} \cdots \sum_{\ell_{4}: \ell_{4}\neq k} \left | \frac{b_{j,\ell_{1}} }{\kappa_{k}-\kappa_{\ell_{1}}} \right | \cdots \left | \frac{b_{j,\ell_{4}}}{\kappa_{k}-\kappa_{\ell_{4}}} \right | \Ep ( | \xi_{\ell_{1}} \cdots \xi_{\ell_{4}} | ) \\
&\quad \lesssim j^{-4\gamma} \sum_{\ell_{1}: \ell_{1} \neq k} \cdots \sum_{\ell_{4}: \ell_{4}\neq k} \frac{\ell_{1}^{-\beta}}{|\kappa_{k}-\kappa_{\ell_{1}}|} \cdots \frac{\ell_{4}^{-\beta}}{|\kappa_{k}-\kappa_{\ell_{4}}|} \Ep ( | \xi_{\ell_{1}} \cdots \xi_{\ell_{4}} | )
\end{align*}
and a repeated application of H\"{o}lder's inequality yields that 
\[
\Ep( | \xi_{\ell_{1}} \cdots \xi_{\ell_{4}} | ) \leq \{\Ep( \xi_{\ell_{1}}^{4})\}^{1/4} \cdots  \{\Ep( \xi_{\ell_{4}}^{4})\}^{1/4} \lesssim \ell_{1}^{-\alpha/2} \cdots \ell_{4}^{-\alpha/2}.
\]
Hence the first term on the right hand side of (\ref{eq: T1}) is 
\[
\lesssim n^{-1} j^{-2\gamma} k^{-\alpha}   \left ( \sum_{\ell: \ell \neq k} \frac{\ell^{-\beta-\alpha/2}}{|\kappa_{k} - \kappa_{\ell}|} \right)^{2} \lesssim n^{-1} j^{-2\gamma} k^{-\alpha},
\]
where the last inequality follows from a similar estimate to (\ref{eq: delicate bound}) together with the assumption that $\beta > \alpha/2+1$. Using a similar argument to bound the second term on the right hand side of (\ref{eq: T1}), we conclude that $\Ep(T_{j,k,1}^{2}) \lesssim n^{-1} j^{-2\gamma} k^{-\alpha}$.

Finally, we shall bound $|T_{j,k,2}|$. 
To this end, observe that, on the event $A_{n}$, 
\[
|T_{j,k,2}| \lesssim  j^{-\gamma}  ||| \hat{K} - K |||  \sum_{\ell: \ell \neq k} \frac{\ell^{-\beta} \hat{v}_{k,\ell}}{|\kappa_{k} - \kappa_{\ell}|^{2}},
\]
where 
\[
\hat{v}_{k,\ell} =\left | \frac{1}{n} \sum_{i=1}^{n} \xi_{i,k} \xi_{i,\ell} - \overline{\xi}_{k} \overline{\xi}_{\ell} \right|.
\]
Then we have
\begin{align*}
\Ep\left \{ \left (\sum_{\ell: \ell \neq k} \frac{\ell^{-\beta}}{|\kappa_{k}-\kappa_{\ell}|^{2}}\hat{v}_{k,\ell} \right )^{2} \right \} 
&\leq \left [\sum_{\ell: \ell \neq k} \frac{\ell^{-\beta}}{|\kappa_{k}-\kappa_{\ell}|^{2}} \{ \Ep(\hat{v}_{k,\ell}^{2})\}^{1/2} \right ]^{2} \\
&\lesssim n^{-1} k^{-\alpha} \left ( \sum_{\ell: \ell \neq k} \frac{\ell^{-\beta-\alpha/2}}{|\kappa_{k}-\kappa_{\ell}|^{2}}  \right )^{2},
\end{align*}
and the far right hand side is $\lesssim n^{-1} (k^{-\alpha} + k^{2\alpha - 2\beta + 4})$, because
\begin{align*}
&\sum_{\ell: \ell \neq k} \frac{\ell^{-\beta-\alpha/2}}{|\kappa_{k}-\kappa_{\ell}|^{2}} = \left ( \sum_{\ell =1}^{[k/C]} + \sum_{\ell = [k/C]+1, \neq k}^{[Ck]} + \sum_{\ell=[Ck]+1}^{\infty} \right ) \frac{\ell^{-\beta-\alpha/2}}{|\kappa_{k}-\kappa_{\ell}|^{2}} \\
&\quad \lesssim \sum_{\ell=1}^{[k/C]} \ell^{3\alpha/2 - \beta} + k^{2\alpha+2} \sum_{\ell = [k/C]+1,\neq k}^{[Ck]} \frac{\ell^{-\beta-\alpha/2}}{|k-\ell|^{2}} + k^{2\alpha} \sum_{\ell=[Ck]+1}^{\infty} \ell^{-\beta-\alpha/2} \\
&\quad \lesssim  1 + k^{3\alpha/2 - \beta +1} \log k +  k^{3\alpha/2 - \beta + 2} + k^{3\alpha/2-\beta+1}  \lesssim 1 + k^{3\alpha/2 - \beta + 2}. 
\end{align*}

Summarizing, using (\ref{eq: eigenvalue bound}) and (\ref{eq: eigenfunction bound}), we have 
$
\sum_{j=1}^{\infty} \sum_{k=1}^{m_{n}} (T_{j,k,1}^{2}+\cdots+T_{j,k,4}^{2}) 
=O_{\Prr} [ n^{-1} + n^{-2}\{ m_{n}^{3} +  m_{n}^{2\alpha - 2\beta + 5} (\log m_{n})^{2} \}   ]$.
Since $m_{n} = o\{n^{1/(2\alpha+2)}\}$, $m_{n}^{3} = o(n)$, so that the the last expression is $O_{\Prr} \{n^{-1}+n^{-2} m_{n}^{2\alpha - 2\beta + 5} (\log m_{n})^{2} \}$.
Furthermore, $m_{n}^{2\alpha - 2\beta + 5} (\log m_{n})^{2}=o (m_{n}^{\alpha+3})= o(n)$ since $\beta > \alpha/2+1$. Hence we conclude that 
\[
\sum_{j=1}^{\infty} \sum_{k=1}^{m_{n}} (T_{j,k,1}^{2}+\cdots+T_{j,k,4}^{2}) = O_{\Prr}(n^{-1}).
\]

\textbf{Step 2}. We wish to bound  $\sum_{j=1}^{\infty} \sum_{k=1}^{m_{n}} \hat{\kappa}_{k}^{-2} (n^{-1} \sum_{i=1}^{n} \varepsilon_{i,j}\hat{\xi}_{i,k})^{2}$. 
Observe that for $1 \leq k \leq m_{n}$, $| \hat{\kappa}_{k}/\kappa_{k} - 1 | \lesssim k^{\alpha} | \hat{\kappa}_{k} - \kappa_{k} | \leq m_{n}^{\alpha} ||| \hat{K} - K ||| = o_{\Prr}(1)$,
from which we have $\max_{1 \leq k \leq m_{n}} | \kappa_{k}/\hat{\kappa}_{k}- 1 | = o_{\Prr}(1)$.
Since conditionally on $X_{1}^{n} = \{ X_{1},\dots,X_{n} \}$, $\varepsilon_{1,j},\dots,\varepsilon_{n,j}$ are independent with mean zero, we have $\Ep  \{ (n^{-1} \sum_{i=1}^{n} \varepsilon_{i,j} \hat{\xi}_{i,k} )^{2} \mid  X_{1}^{n} \} = n^{-2} \sum_{i=1}^{n} \Ep(\varepsilon_{i,j}^{2} \mid X_{1}^{n}) \hat{\xi}_{i,k}^{2}$
Further, since by the monotone convergence theorem for conditional expectation and Bessel's inequality, 
$\sum_{j=1}^{\infty} \Ep (\varepsilon_{i,j}^{2} \mid X_{1}^{n}) =\Ep ( \sum_{j=1}^{\infty} \varepsilon_{i,j}^{2} \mid  X_{1}^{n} ) \leq \Ep ( \| \mE_{i} \|^{2} \mid X_{1}^{n} ) = \Ep( \| \mE_{i} \|^{2} \mid X_{i}) \leq C_{1}$,
we have 
\begin{align*}
&\Ep  \left \{ \sum_{j=1}^{\infty} \sum_{k=1}^{m_{n}} \hat{\kappa}_{k}^{-2} \left  (\frac{1}{n} \sum_{i=1}^{n} \varepsilon_{i,j} \hat{\xi}_{i,k} \right )^{2} \ \Bigg | \ X_{1}^{n} \right \}   \lesssim n^{-1} \sum_{k=1}^{m_{n}} \hat{\kappa}_{k}^{-1} \\
&\quad \leq n^{-1} \{ 1+o_{\Prr}(1) \} \sum_{k=1}^{m_{n}} \kappa_{k}^{-1} = O_{\Prr} (n^{-1} m_{n}^{\alpha+1}).
\end{align*}
This yields that 
\[
\sum_{j=1}^{\infty} \sum_{k=1}^{m_{n}} \hat{\kappa}_{k}^{-2}\left  (\frac{1}{n}\sum_{i=1}^{n} \varepsilon_{i,j} \hat{\xi}_{i,k} \right)^{2} = O_{\Prr} (n^{-1}m_{n}^{\alpha+1}).
\]
Summarizing, we conclude that 
\[
\sum_{j=1}^{\infty} \sum_{k=1}^{m_{n}} (\hat{b}_{j,k} - b_{j,k})^{2} = O_{\Prr} (n^{-1} m_{n}^{\alpha+1} ).
\]

\textbf{Step 3}. Conclusion. Recall that $\hat{b} = \sum_{j=1}^{\infty} \sum_{k=1}^{m_{n}} \hat{b}_{j,k} (\phi_{j} \otimes \hat{\phi}_{k})$, and observe that 
\[
\hat{b} - b = \sum_{j=1}^{\infty} \sum_{k=1}^{m_{n,2}} (\hat{b}_{j,k} - b_{j,k}) (\phi_{j} \otimes \hat{\phi}_{k}) + \sum_{j=1}^{\infty} \sum_{k=1}^{m_{n}} b_{j,k} \{ \phi_{j} \otimes ( \hat{\phi}_{k} - \phi_{k}) \}+ B_{n},
\]
where $B_{n} = b - \sum_{j=1}^{\infty} \sum_{k=1}^{m_{n}} b_{j,k} (\phi_{j} \otimes \phi_{k})$. Since $||| B_{n} |||^{2} = \sum_{j=1}^{\infty} \sum_{k> m_{n}}b_{j,k}^{2} = O(m_{n}^{-2\beta+1})$, 
\[
||| \hat{b} - b |||^{2} =
 O_{\Prr} \left ( n^{-1} m_{n}^{\alpha+1} +  m_{n}^{-2\beta + 1} \right ) + \iint \left [  \sum_{j=1}^{\infty} \sum_{k=1}^{m_{n}} b_{j,k} \{ \phi_{j} \otimes ( \hat{\phi}_{k} - \phi_{k}) \} \right ]^{2}.
\]
Now, observe that, using Parseval's identity, the second term on the right hand side is
\begin{align*}
 \sum_{j=1}^{\infty} \int \left \{ \sum_{k=1}^{m_{n}}  b_{j,k} ( \hat{\phi}_{k} - \phi_{k}) \right \}^{2} 
&\leq m_{n} \sum_{j=1}^{\infty} \sum_{k=1}^{m_{n}} b_{j,k}^{2} \| \hat{\phi}_{k} - \phi_{k} \|^{2} \\
&\lesssim m_{n} \sum_{k=1}^{m_{n}} k^{-2\beta} \| \hat{\phi}_{k} - \phi_{k} \|^{2}, 
\end{align*}
which is $O_{\Prr}(n^{-1}m_{n})$.
This completes the proof for the first assertion. The second assertion follows directly from the first assertion.  \qed

\subsection{Proof of Theorem \ref{thm: main2}}

The proof is parallel to that of Theorem \ref{thm: main1}. We freely use the results in the proof of Theorem \ref{thm: main1}. Since $\tilde{b}$ is invariant with respect to choices of signs of $\hat{\phi}_{k}$'s, it is without loss of generality to assume (\ref{eq: sign of phi}). 
Let $\overline{m}_{n} = \max \{ m_{n,1},m_{n,2} \}$, and define the event 
\[
A_{n}' =  \{ |\hat{\kappa}_{k}-\kappa_{\ell}| \geq  |\kappa_{k}-\kappa_{\ell}|/2 \ \text{for all} \ 1 \leq k \leq \overline{m}_{n} \ \text{and for all} \ \ell \neq k \},
\] 
for which we have $\Prr (A_{n}') \to 0$ since $\overline{m}_{n} = o\{n^{1/(2\alpha+2)}\}$. 

Expand $\mE_{i} = Y_{i}-\Ep(Y_{i} \mid X_{i})$ as $\mE_{i} = \sum_{j} \hat{\varepsilon}_{i,j} \hat{\phi}_{j}$ with $\hat{\varepsilon}_{i,j} = \int \mE_{i} \hat{\phi}_{j}$.
Let $\hat{\eta}_{i,j}^{c} = \hat{\eta}_{i,j} - n^{-1}\sum_{i'=1}^{n} \hat{\eta}_{i',j}$ and $\hat{\varepsilon}_{i,\ell}^{c} = \hat{\varepsilon}_{i,\ell} - n^{-1} \sum_{i'=1}^{n} \hat{\varepsilon}_{i',\ell}$; observe that 
\[
\hat{\eta}_{i,j}^{c}=\sum_{k=1}^{\infty} \breve{b}^{*}_{j,k} \hat{\xi}_{i,k} + \hat{\varepsilon}_{i,j}^{c},
\]
where $\breve{b}_{j,k}^{*} = \iint b (\hat{\phi}_{j} \otimes \hat{\phi}_{k})$. Hence, using the relation in (\ref{eq: orthogonality}), we have $\hat{\kappa}_{k} \tilde{b}_{j,k} = \hat{\kappa}_{k} \breve{b}^{*}_{j,k} + n^{-1} \sum_{i=1}^{n} \hat{\varepsilon}_{i,j} \hat{\xi}_{i,k}$, 
which yields that $(\tilde{b}_{j,k}-b_{j,k})^{2} \lesssim (\breve{b}^{*}_{j,k}-b_{j,k})^{2} +  \hat{\kappa}_{k}^{-2}   (n^{-1}\sum_{i=1}^{n} \hat{\varepsilon}_{i,j} \hat{\xi}_{i,k} )^{2}$.

Observe that 
\begin{align*}
\breve{b}^{*}_{j,k}-b_{j,k} &= \iint b (\hat{\phi}_{j} \otimes \hat{\phi}_{k} - \phi_{j} \otimes \phi_{k} )\\
&= \iint b \{ (\hat{\phi}_{j} - \phi_{j}) \otimes \phi_{k} \} + \iint  b \{ \phi_{j} \otimes (\hat{\phi}_{k} - \phi_{k}) \}  + \iint b \{ (\hat{\phi}_{j}-\phi_{j}) \otimes (\hat{\phi}_{k}-\phi_{k}) \} \\
&=:I_{j,k} + II_{j,k} + III_{j,k}. 
\end{align*}
Step 1 in the proof of Theorem \ref{thm: main1} shows that $\sum_{j=1}^{m_{n,1}} \sum_{k=1}^{m_{n,2}} II_{j,k}^{2} = O_{\Prr}(n^{-1})$, and likewise we have  $\sum_{j=1}^{m_{n,1}} \sum_{k=1}^{m_{n,2}} I_{j,k}^{2} = O_{\Prr}(n^{-1})$. Furthermore, using (\ref{eq: difference}), observe that $
\iint b \{ (\hat{\phi}_{j} - \phi_{j}) \otimes (\hat{\phi}_{k}-\phi_{k}) \} = \sum_{p,q} b_{p,q} \hat{w}_{j,p} \hat{w}_{k,q}$,
where 
\[
\hat{w}_{j,p} = 
\begin{cases}
(\hat{\kappa}_{j} - \kappa_{p})^{-1} \iint (\hat{K} - K) (\hat{\phi}_{j} \otimes \phi_{p})  & \text{if} \ p \neq j \\
\int (\hat{\phi}_{j} -\phi_{j}) \phi_{j} & \text{if} \ p=j
\end{cases}
.
\]
For each $p \neq j$, on the event $A_{n}'$, $|\hat{w}_{j,p}| \lesssim |\kappa_{j} - \kappa_{p}|^{-1} ||| \hat{K} - K |||$,
which yields that on the event $A_{n}'$, 
\begin{align*}
\left | \sum_{p,q}b_{p,q} \hat{w}_{j,p} \hat{w}_{k,q} \right | &\lesssim \left ( ||| \hat{K} - K ||| \sum_{p: p\neq j}\frac{p^{-\gamma}}{|\kappa_{j} - \kappa_{p}|} + j^{-\gamma} \| \hat{\phi}_{j} - \phi_{j} \| \right ) \\
&\quad \times  \left ( ||| \hat{K} - K ||| \sum_{q:q \neq k} \frac{q^{-\beta}}{|\kappa_{k} - \kappa_{q}|} + k^{-\beta} \| \hat{\phi}_{k} - \phi_{k} \| \right ) \\
&\lesssim \{ ||| \hat{K} - K ||| ( 1 +  j^{\alpha-\gamma+1} \log j ) +  j^{-\gamma} \| \hat{\phi}_{j} - \phi_{j} \| \} \\
&\quad \times  \{ ||| \hat{K} - K ||| (1+ k^{\alpha-\beta+1} \log k )  +  k^{-\beta} \| \hat{\phi}_{k} - \phi_{k} \| \}.
\end{align*}
Therefore, we have
\begin{align*}
&\sum_{j=1}^{m_{n,1}} \sum_{k=1}^{m_{n,2}} \left ( \sum_{p,q}b_{p,q} \hat{w}_{j,p} \hat{w}_{k,q} \right )^{2} \\
&= O_{\Prr} \left [  n^{-2} \{  m_{n,1}+ m_{n,1} ^{2\alpha-2\gamma+3} (\log m_{n,1})^{2} \}  \{  m_{n,2} + m_{n,2} ^{2\alpha-2\beta+3} (\log m_{n,2})^{2} \}  \right ].
\end{align*}
Since $\beta > \alpha/2+1$ and $\gamma > \alpha/2+1$, the last expression is $o_{\Prr} (n^{-2} m_{n,1}^{\alpha+1}m_{n,2}^{\alpha+1}) = o_{\Prr}(n^{-1})$. 
So we conclude that $\sum_{j=1}^{m_{n,1}} \sum_{k=1}^{m_{n,2}} (\breve{b}^{*}_{j,k} - b_{j,k})^{2} =O_{\Prr} (n^{-1})$. 

Next, since conditionally on $X_{1}^{n} =\{ X_{1},\dots,X_{n} \}$, $\hat{\varepsilon}_{1,j},\dots,\hat{\varepsilon}_{n,j}$ are independent with mean zero, we have 
\[
\Ep \left \{ \left ( \frac{1}{n}\sum_{i=1}^{n} \hat{\varepsilon}_{i,j} \hat{\xi}_{i,k} \right )^{2} \ \Bigg| \ X_{1}^{n} \right \} = \frac{1}{n^{2}} \sum_{i=1}^{n} \Ep(\hat{\varepsilon}_{i,j}^{2} \mid X_{1}^{n}) \hat{\xi}_{i,k}^{2}.
\] 
Further, since by Bessel's inequality, 
$\sum_{j=1}^{m_{n,1}} \Ep (\hat{\varepsilon}_{i,j}^{2} \mid X_{1}^{n}) \leq \Ep ( \sum_{j=1}^{m_{n,1}} \hat{\varepsilon}_{i,j}^{2} \mid X_{1}^{n} ) \leq \Ep ( \| \mE_{i} \|^{2} \mid X_{1}^{n} ) = \Ep (\| \mE_{i} \|^{2} \mid X_{i}) \leq C_{1}$, 
we have, using the fact that $\max_{1 \leq k \leq m_{n,2}} | \kappa_{k}/\hat{\kappa}_{k}- 1 | = o_{\Prr}(1)$, 
\begin{align*}
&\Ep \left \{ \sum_{j=1}^{m_{n,1}} \sum_{k=1}^{m_{n,2}} \frac{1}{\hat{\kappa}_{k}^{2}} \left (\frac{1}{n} \sum_{i=1}^{n} \hat{\varepsilon}_{i,j} \hat{\xi}_{i,k} \right)^{2} \ \Bigg| \ X_{1}^{n} \right \}  \lesssim n^{-1} \sum_{k=1}^{m_{n,2}} \hat{\kappa}_{k}^{-1} \\
&\quad \leq n^{-1} \{ 1+o_{\Prr}(1) \} \sum_{k=1}^{m_{n,2}} \kappa_{k}^{-1} = O_{\Prr} (n^{-1} m_{n,2}^{\alpha+1}). 
\end{align*}
This yields that $\sum_{j=1}^{m_{n,1}} \sum_{k=1}^{m_{n,2}} \hat{\kappa}_{k}^{-2} (n^{-1} \sum_{i=1}^{n} \hat{\varepsilon}_{i,j} \hat{\xi}_{i,k} )^{2} = O_{\Prr} (n^{-1}m_{n,2}^{\alpha+1})$. 

Summarizing, we conclude that $\sum_{j=1}^{m_{n,1}} \sum_{k=1}^{m_{n,2}} (\tilde{b}_{j,k} - b_{j,k})^{2} = O_{\Prr} (n^{-1} m_{n,2}^{\alpha+1} )$.

Recall that $\tilde{b} = \sum_{j=1}^{m_{n,1}} \sum_{k=1}^{m_{n,2}} \tilde{b}_{j,k} (\hat{\phi}_{j} \otimes \hat{\phi}_{k})$, and observe that 
\[
\tilde{b} - b = \sum_{j=1}^{m_{n,1}} \sum_{k=1}^{m_{n,2}} (\tilde{b}_{j,k} - b_{j,k}) (\hat{\phi}_{j} \otimes \hat{\phi}_{k}) + \sum_{j=1}^{m_{n,1}} \sum_{k=1}^{m_{n,2}} b_{j,k} (\hat{\phi}_{j} \otimes \hat{\phi}_{k} - \phi_{j} \otimes \phi_{k}) + B'_{n},
\]
where $B'_{n} = b - \sum_{j=1}^{m_{n,1}} \sum_{k=1}^{m_{n,2}} b_{j,k} (\phi_{j} \otimes \phi_{k})$. So 
\begin{align*}
||| \tilde{b} - b |||^{2} &\lesssim  \sum_{j=1}^{m_{n,1}} \sum_{k=1}^{m_{n,2}}  (\tilde{b}_{j,k} - b_{j,k})^{2} + \iint \left \{ \sum_{j=1}^{m_{n,1}} \sum_{k=1}^{m_{n,2}} b_{j,k} (\hat{\phi}_{j} \otimes \hat{\phi}_{k} - \phi_{j} \otimes \phi_{k}) \right \}^{2} \\
&\quad + \left ( \sum_{j > m_{n,1}} \sum_{k=1}^{m_{n,2}} + \sum_{j=1}^{m_{n,1}} \sum_{k> m_{n,2}} + \sum_{j > m_{n,1}} \sum_{k > m_{n,2}} \right ) b_{j,k}^{2} \\
&= O_{\Prr} \left ( n^{-1} m_{n,2}^{\alpha+1} + m_{n,1}^{-2\gamma+1} + m_{n,2}^{-2\beta + 1} \right ) \\
&\quad + \iint \left \{ \sum_{j=1}^{m_{n,1}} \sum_{k=1}^{m_{n,2}} b_{j,k} (\hat{\phi}_{j} \otimes \hat{\phi}_{k} - \phi_{j} \otimes \phi_{k}) \right \}^{2}.
\end{align*}
Using the decomposition 
\[
\hat{\phi}_{j} \otimes \hat{\phi}_{k} - \phi_{j} \otimes \phi_{k} = (\hat{\phi}_{j} - \phi_{j} )\otimes \phi_{k} + \phi_{j} \otimes  (\hat{\phi}_{k} - \phi_{k})+  (\hat{\phi}_{j}-\phi_{j}) \otimes (\hat{\phi}_{k}-\phi_{k}),
\]
we have 
\begin{align*}
&\sum_{j=1}^{m_{n,1}} \sum_{k=1}^{m_{n,2}} b_{j,k} (\hat{\phi}_{j} \otimes \hat{\phi}_{k} - \phi_{j} \otimes \phi_{k}) \\
&=\sum_{j=1}^{m_{n,1}} (\hat{\phi}_{j} - \phi_{j} ) \otimes \left (\sum_{k=1}^{m_{n,2}} b_{j,k} \phi_{k} \right )  +\sum_{k=1}^{m_{n,2}} \left (\sum_{j=1}^{m_{n,1}} b_{j,k} \phi_{j} \right ) \otimes (\hat{\phi}_{k} - \phi_{k} )  \\
&\quad + \sum_{j=1}^{m_{n,1}} \sum_{k=1}^{m_{n,2}} b_{j,k} \{ (\hat{\phi}_{j}-\phi_{j}) \otimes (\hat{\phi}_{k}-\phi_{k}) \}. 
\end{align*}
Observe that 
\begin{align*}
&\iint \left \{ \sum_{j=1}^{m_{n,1}} (\hat{\phi}_{j} - \phi_{j} ) \otimes \left (\sum_{k=1}^{m_{n,2}} b_{j,k} \phi_{k} \right ) \right \}^{2} \leq m_{n,1} \sum_{j=1}^{m_{n,1}} \| \hat{\phi}_{j} - \phi_{j} \|^{2} \sum_{k=1}^{m_{n,2}} b_{j,k}^{2} \\
&\quad \lesssim m_{n,1} \sum_{j=1}^{m_{n,1}} j^{-2\gamma} \| \hat{\phi}_{j} - \phi_{j} \|^{2} = O_{\Prr} (n^{-1}m_{n,1}). 
\end{align*}
Likewise, we have 
\[
\iint \left \{ \sum_{k=1}^{m_{n,2}} \left (\sum_{j=1}^{m_{n,1}} b_{j,k} \phi_{j} \right ) \otimes (\hat{\phi}_{k} - \phi_{k} ) \right \}^{2} = O_{\Prr} (n^{-1}m_{n,2}).
\]
Finally, we have 
\begin{align*}
&\iint \left [ \sum_{j=1}^{m_{n,1}} \sum_{k=1}^{m_{n,2}} b_{j,k} \{ (\hat{\phi}_{j}-\phi_{j}) \otimes (\hat{\phi}_{k}-\phi_{k}) \} \right ] \\
&\leq m_{n,1} m_{n,2} \sum_{j=1}^{m_{n,1}} \sum_{k=1}^{m_{n,2}} b_{j,k}^{2} \| \hat{\phi}_{j} - \phi_{j} \|^{2} \| \hat{\phi}_{k} - \phi_{k} \|^{2} = O_{\Prr} ( n^{-2} m_{n,1}m_{n,2} ). 
\end{align*}

Therefore, we conclude that
\[
||| \tilde{b} - b |||^{2} = O_{\Prr} \left \{ n^{-1} (m_{n,1} + m_{n,2}^{\alpha+1}) +  m_{n,1}^{-2\gamma+1} + m_{n,2}^{-2\beta + 1} \right \}.
\]
The second assertion follows directly from the first assertion. This completes the proof. 

\subsection{Proof of Theorem \ref{thm: lower bound}}

The proof is inspired by that of (3.6) in \cite{HaHo07}; the current proof relies on Assouad's lemma \citep[][Lemma 2.12]{Ts03} and  Theorem 2.12 in \cite{Ts03}. To apply those results, we have to construct a sequence of conditional distributions of $Y$ given $X$, to which end we employ the theory of Gaussian measures on Banach spaces; see, e.g., \cite{St11}, Chapter VIII. 

For any $b \in L^{2}(I^{2})$ and $x \in L^{2}(I)$, let $P_{b,x}$ denote the distribution of $\int_{I} b(\cdot, t) x(t)dt + \mE(\cdot)$, and let $P_{0}$ denote the distribution of $\mE$. Those distributions are defined on the Borel $\sigma$-field of $L^{2}(I)$.
Associated to $\mE$, the Cameron-Martin space is given by 
\[
H = \left \{ h = \sum_{j} h_{j} \phi_{j} : \sum_{j} \frac{h_{j}^{2}}{\lambda_{j}} < \infty \right \}
\]
equipped with the inner product 
\[
\langle h,g \rangle_{H} = \sum_{j} \frac{h_{j}g_{j}}{\lambda_{j}}, \ h=\sum_{j} h_{j} \phi_{j}, g=\sum_{j}g_{j} \phi_{j} \in H. 
\]
Let $b = \sum_{j,k} b_{j,k} \phi_{j}\otimes \phi_{k}$ and $x = \sum_{k} x_{k}\phi_{k}$; then  $P_{b,x}$ is absolutely continuous with respect to $P_{0}$ if and only if
\[
\left \| \int_{I} b(\cdot,t)x(t) dt \right \|_{H}^{2} = \sum_{j} \frac{(\sum_{k} b_{j,k} x_{k})^{2}}{\lambda_{j}} < \infty,
\]
and its Radon-Nikodym derivative is given by the Cameron-Martin formula
\[
p_{b,x}(y) = \frac{d P_{b,x}}{dP_{0}}(y)  
= \exp \left \{ -\sum_{j} \frac{(\sum_{k} b_{j,k} x_{k})^{2}}{2\lambda_{j}} + \sum_{j} \frac{y_{j} \sum_{k} b_{j,k} x_{k}}{\lambda_{j}} \right \}, 
\]
where $y=\sum_{j} y_{j} \phi_{j}$. See Theorem 8.2.9 in \cite{St11}. Denote by $Q$ the distribution of $X$; then the joint distribution of $(X,Y)$ is given by  $p_{b,x} (y) dP_{0}(y) dQ(x)$. 

Now, let $\nu_{n} = [n^{1/(\alpha+2\beta)}]$, and 
\[
b^{\theta} = \sum_{k=\nu_{n}+1}^{2\nu_{n}} k^{-\beta} \theta_{k-\nu_{n}} (\phi_{1} \otimes \phi_{k}),
\]
where $\theta = (\theta_{1},\dots,\theta_{\nu_{n}}) \in \{ 0,1 \}^{\nu_{n}}$. Then $b^{\theta} \in \mB (\beta,\gamma,C_{1})$ and $b^{\theta}_{j,k} = 0$ for all $j \geq 2$, so that 
\begin{align*}
p_{b^{\theta},x}(y)  &= \exp \left \{ - \frac{(\sum_{k=\nu_{n}+1}^{2\nu_{n}} k^{-\beta} \theta_{k-\nu_{n}} x_{k})^{2}}{2\lambda_{1}} + \frac{y_{1} \sum_{k=\nu_{n}+1}^{2\nu_{n}} k^{-\beta} \theta_{k-\nu_{n}}x_{k}}{\lambda_{1}} \right\}.
\end{align*}
Define $\tilde{p}_{\theta,x}(y) = p_{b^{\theta},x}(y)$ and $d\tilde{P}_{\theta} (x,y)= \tilde{p}_{\theta,x} (y) dP_{0}(y)dQ(x)$ for each $\theta = (\theta_{1},\dots,\theta_{\nu_{n}}) \in \{ 0,1 \}^{\nu_{n}}$, and let $(X_{1},Y_{1}),\dots,(X_{n},Y_{n})$ be i.i.d. from $\tilde{P}_{\theta}$. 

For any estimator $\overline{b}^{n} = \sum_{j,k} \overline{b}^{n}_{j,k} (\phi_{j} \otimes \phi_{k})$ of $b^{\theta}$, we have by Bessel's inequality, 
\begin{align*}
||| \overline{b}^{n} - b^{\theta} |||^{2} &\geq \sum_{k=\nu_{n}+1}^{2\nu_{n}} (\overline{b}^{n}_{1,k} - k^{-\beta}\theta_{k-\nu_{n}})^{2} \\
&\geq \frac{1}{4} \sum_{k=\nu_{n}+1}^{2\nu_{n}} k^{-2\beta} (\overline{\theta}^{n}_{k-\nu_{n}}-\theta_{k-\nu_{n}})^{2} \\
&\geq \frac{(2\nu_{n})^{-2\beta}}{4} \sum_{k=1}^{\nu_{n}} | \overline{\theta}^{n}_{k} - \theta_{k} |,
\end{align*}
where 
\[
\overline{\theta}^{n}_{k-\nu_{n}} = \arg \min_{\vartheta \in \{ 0,1 \}} (k^{\beta} \overline{b}^{n}_{1,k} - \vartheta)^{2} 
=
\begin{cases}
0 & \text{if} \ k^{\beta} \overline{b}^{n}_{1,k} \leq 1/2 \\
1 & \text{if} \ k^{\beta} \overline{b}^{n}_{1,k} > 1/2
\end{cases}
.
\]
Indeed, since $(k^{\beta} \overline{b}^{n}_{1,k} - \overline{\theta}^{n}_{k-\nu_{n}})^{2} \leq (k^{\beta}\overline{b}^{n}_{1,k} - \theta_{k-\nu_{n}})^{2}$ by the definition of $\overline{\theta}^{n}_{k-\nu_{n}}$,  
\begin{align*}
(\overline{\theta}^{n}_{k-\nu_{n}}-\theta_{k-\nu_{n}})^{2} &\leq 2(k^{\beta} \overline{b}^{n}_{1,k} - \overline{\theta}^{n}_{k-\nu_{n}})^{2} + 2(k^{\beta}\overline{b}^{n}_{1,k} - \theta_{k-\nu_{n}})^{2} \\
&\leq 4(k^{\beta}\overline{b}^{n}_{1,k} - \theta_{k-\nu_{n}})^{2}.
\end{align*}

For any $\theta, \theta' \in \{ 0,1 \}^{\nu_{n}}$, let $\rho (\theta, \theta') = \sum_{k=1}^{\nu_{n}} | \theta_{k} - \theta'_{k} |$ denote the Hamming distance. Then we have 
\[
\Prr_{\theta} \left \{  ||| \overline{b}^{n} - b^{\theta} |||^{2} \geq \frac{(2\nu_{n})^{-2\beta}}{4}c \right \} \geq  \Prr_{\theta} \left \{  \rho (\overline{\theta}^{n}, \theta) \geq c \right \}
\]
for any $\theta \in \{ 0,1 \}^{\nu_{n}}$ and any constant $c>0$, where $\Prr_{\theta}$ denotes the probability under $\theta$. 
To lower bound the right hand side, we calculate the Kullback-Leibler divergence
\[
K(\tilde{P}_{\theta},\tilde{P}_{\theta'}) = \int \log \frac{d\tilde{P}_{\theta}}{d\tilde{P}_{\theta'}} d\tilde{P}_{\theta}
\]
for any $\theta, \theta' \in \{ 0,1 \}^{\nu_{n}}$ with $\rho (\theta,\theta') = 1$. 
Suppose that $\theta_{k} \neq \theta'_{k}$ for some $1 \leq k \leq \nu_{n}$ and $\theta_{\ell} = \theta'_{\ell}$ for all $\ell \neq k$. Then a straightforward calculation shows that 
\[
K(\tilde{P}_{\theta},\tilde{P}_{\theta'}) = \Ep_{\theta} \left \{ \log \frac{\tilde{p}_{\theta,X}(Y)}{\tilde{p}_{\theta',X}(Y)} \right \} =\frac{(\nu_{n}+k)^{-\alpha-2\beta}}{2\lambda_{1}} \leq \frac{(\nu_{n}+1)^{-\alpha-2\beta}}{2\lambda_{1}} \leq \frac{1}{2\lambda_{1}n},
\]
which yields that
\[
K(\tilde{P}^{\otimes n}_{\theta},\tilde{P}^{\otimes n}_{\theta'}) = nK(\tilde{P}_{\theta},\tilde{P}_{\theta'}) \leq \frac{1}{2\lambda_{1}}.
\]
Now, applying Assouad's lemma and Theorem 2.12 in \cite{Ts03}, we have 
\[
 \max_{\theta \in \{ 0,1 \}^{\nu_{n}} }\Ep_{\theta} \{\rho (\overline{\theta}^{n}, \theta) \} \geq \frac{\nu_{n}}{4} e^{-1/(2\lambda_{1})},
\]
where $\Ep_{\theta}$ denotes the expectation under $\theta$. Choose $\theta \in \{ 0,1 \}^{\nu_{n}}$ at which the maximum on the left hand side is attained, and observe that $\rho (\overline{\theta}^{n},\theta) \leq \nu_{n}$. The Paley-Zygmund inequality then yields that
\begin{align*}
\Prr_{\theta} \left \{ \rho (\overline{\theta}^{n}, \theta)  \geq \frac{\nu_{n}}{8} e^{-1/(2\lambda_{1})} \right \} 
&\geq 
\Prr_{\theta} \left [ \rho (\overline{\theta}^{n},\theta) \geq \frac{1}{2} \Ep_{\theta}\{\rho(\overline{\theta}^{n},\theta)\} \right ] \\
&\geq \frac{1}{4} \frac{[\Ep_{\theta}\{\rho (\overline{\theta}^{n},\theta)\}]^{2}}{\Ep\{\rho (\overline{\theta}^{n},\theta)^{2}\}} \\
&\geq \frac{1}{16} e^{-1/(2\lambda_{1})}.
\end{align*}
Therefore 
\[
\max_{\theta \in \{0,1\}^{\nu_{n}}} \Prr_{\theta} \left \{  ||| \overline{b}^{n} - b^{\theta} |||^{2} \geq \frac{\nu_{n}^{-2\beta+1}}{2^{2\beta+5}}e^{-1/(2\lambda_{1})} \right \} \geq \frac{1}{16} e^{-1/(2\lambda_{1})}.
\]
Since $\nu_{n}^{-2\beta+1} \sim n^{-(2\beta-1)/(\alpha+2\beta)}$, the proof is completed. \qed

\end{document}